\documentclass[11pt]{amsart}
\usepackage[margin=1.1in]{geometry}
\usepackage{amsmath,amssymb,amsthm,booktabs}
\usepackage[hidelinks]{hyperref}

\newtheorem{theorem}{Theorem}[section]
\newtheorem{lemma}[theorem]{Lemma}
\newtheorem{proposition}[theorem]{Proposition}
\newtheorem{corollary}[theorem]{Corollary}
\theoremstyle{definition}
\newtheorem{definition}[theorem]{Definition}
\theoremstyle{remark}
\newtheorem{remark}[theorem]{Remark}

\numberwithin{equation}{section}
\newcommand{\Z}{\mathbb{Z}}
\newcommand{\R}{\mathbb{R}}
\newcommand{\Q}{\mathbb{Q}}
\newcommand{\gap}{g}
\allowdisplaybreaks

\title[Positive values of inhomogeneous forms of type $(1,4)$]{Positive values of inhomogeneous quadratic forms of type $(1,4)$:\\ the constant $\Gamma_{1,4}=8$}
\author{Vineet Kahlon}
\date{\today}

\begin{document}

\begin{abstract}
Let $Q$ be a real indefinite quadratic form of type $(1,4)$ with determinant $D\neq 0$. We prove that for every $c\in\R^5$ there exists $x\in\Z^5$ such that $0<Q(x+c)\le(8|D|)^{1/5}$, and that the constant $8$ is optimal. This establishes a conjecture of Bambah, Dumir and Hans-Gill and, together with earlier work, completes the determination of the sharp constants $\Gamma_{r,s}$ for positive inhomogeneous values of nondegenerate indefinite quadratic forms. The proof combines classical reduction and gap theorems with uniform bounds for circular gaps of quadratic residues. A central step identifies the exceptional parameters in certain resonant reduced families: an explicit ladder of integer witnesses resolves their punctured neighbourhoods, leaving only three rational parameter points, subject to stated coefficient conditions. These are treated by retaining the coupling in a preceding, less reduced form. The remaining covering statements are verified by finite certificates in exact rational arithmetic.
\end{abstract}

\maketitle

\section{Introduction}\label{sec:intro}

Let $Q(x_1,\dots,x_n)$ be a real indefinite quadratic form of type $(r,s)$, $r+s=n$, with determinant $D\neq0$ (the determinant of its symmetric Gram matrix). Following Bambah, Dumir and Hans-Gill, let $\Gamma_{r,s}$ be the infimum of the constants $\Gamma$ with the following property: for every such form and every $c\in\R^n$ there is $x\in\Z^n$ with
\[
0<Q(x+c)\le(\Gamma|D|)^{1/n}.
\]
These constants measure a uniform problem in the geometry of numbers: after normalisation by the determinant, how long must an interval $(0,T]$ be to contain a value of the form on every translate $\Z^n+c$ of the integer lattice? The uniformity is over both the form and the translate, and the requirement of positivity distinguishes the problem from that of finding values of small absolute value, as for the inhomogeneous minimum.

For nondegenerate indefinite forms, every constant $\Gamma_{r,s}$ except one is now known. Work between 1947 and 1984 determined all the constants in dimensions at most four and three of the four indefinite types in dimension five~\cite{BDH84,Barnes61,Blaney48,Blaney50,DH47,Dumir67,Dumir68,Dumir68b,DHG81,HGR80,HGR81}. For $n\ge6$, Bambah, Dumir and Hans-Gill conjectured a formula depending only on $n$ and the signature~\cite{BDH80,BDH00}, and its proof was completed through the results of~\cite{AG88,AG91a,AG91b,BDH83,BDH84,DHS95,Flahive88,SD94}. The signature-modulo-eight theorem of Dumir, Hans-Gill and Woods~\cite{DHW94} (in fixed dimension, $\Gamma_{r,s}$ depends only on the signature modulo~$8$) provides a key organising principle. We refer to the surveys~\cite{BDH00} and~\cite[\S1]{BKR}; the latter records $\Gamma_{1,4}$ as the only case left.

For type $(1,4)$, Bambah, Dumir and Hans-Gill conjectured that $\Gamma_{1,4}=8$. Jackson's work~\cite{Jackson} gives $\Gamma_{1,4}\le32$. Dumir and Sehmi~\cite{DS94} established the lower bound $8$ and an upper bound $16$; Raka and Rani~\cite{RR97} subsequently improved the upper bound to $12$. Bhardwaj, Kathuria and Raka~\cite{BKR} proved the conjectured bound when the second shift coordinate is non-integral, when it is integral and the minimum $a$ of the associated ternary form satisfies $a\ge 1/2$, and in several further subcases with $a<1/2$. In the subcases they left open they obtained weaker estimates: the bound $32/3$ in place of $8$ for three of them, and $8.486$ for a fourth. We prove the sharp bound in all remaining cases. Our argument follows the classical canonical reduction but does not invoke their case theorems. With Theorem~\ref{thm:main}, the determination of the constants $\Gamma_{r,s}$ for nondegenerate indefinite forms, begun by Davenport and Heilbronn in 1947, is complete.

\begin{theorem}\label{thm:main}
Let $Q$ be a real quadratic form of type $(1,4)$ with determinant $D\ne0$. For every $c\in\R^5$ there exists $x\in\Z^5$ with
\[
0<Q(x+c)\le(8|D|)^{1/5}.
\]
The constant $8$ cannot be replaced by any smaller number. Hence $\Gamma_{1,4}=8$.
\end{theorem}

Sharpness is shown by an explicit extremal form and shift (Section~\ref{sec:end}); we do not attempt to classify all cases of equality.

\subsection*{The difficulty}
The problem is not merely to find a small value of $Q$, but a \emph{positive} value below the sharp threshold, uniformly over all translates $x+c$. Classical reduction and one-variable gap theorems settle most parameter ranges. At certain rational resonances, however, the reduced inequalities have genuine exceptional configurations, so repeating the same one-variable argument cannot finish the proof. These \emph{pinch points} are obstructions to a reduced auxiliary problem, not counterexamples to Theorem~\ref{thm:main}: the reduction has discarded information that the original form still carries. The strategy is therefore to determine exactly where the reduced problem fails, and to restore the discarded coupling only there.

\subsection*{The new ingredients}
The classical lemmas are the inputs of the proof; its new part is the analysis of the exceptional configurations they leave. In the resonant families of Section~\ref{sec:22} we determine the exact failure set of the twice-completed form $G$. At the resonances $t=1/4$ and $t=1/6$ it consists of three rational values of the unrestricted linear parameters $(\alpha,\kappa)$, subject to the coefficient conditions of Theorem~\ref{thm:fail}. Near these points the successful witnesses cannot be chosen uniformly bounded, so no finite search on the parameter torus closes their neighbourhoods. Explicit infinite \emph{ladders} of integer witnesses, given by polynomial parametrisations whose consecutive members satisfy simple handoff inequalities, certify punctured neighbourhoods of the exceptional points (Lemma~\ref{lem:ladder}). Away from these certified neighbourhoods, compactness yields finite dyadic witness covers (Proposition~\ref{prop:complete}). At the three points themselves we return to the once-completed form $F$, keep the coupling between its coefficients, and conclude with three small Farkas trees.

These arguments are formulated independently of the type $(1,4)$. For compact parameter families with strict witness coverage away from certified ladder neighbourhoods, the completeness criterion gives a finite dyadic cover; the ladder lemma supplies such neighbourhoods when its rung and handoff conditions can be established. Uniform residue-gap bounds (Proposition~\ref{prop:gapexact}, Lemmas~\ref{lem:twopiece} and~\ref{lem:highq}) and a dimension-independent quotient-pruning criterion (Lemma~\ref{lem:prune}) organise the remaining reductions.

Table~\ref{tab:status} compares, case by case, the constants obtained in~\cite{BKR} with those obtained here; the case labels are explained in the outline below.

\begin{table}[ht]
\centering\small
\begin{tabular}{lll}
\toprule
Case (after the canonical reduction) & Bhardwaj--Kathuria--Raka~\cite{BKR} & This paper: $8$\\
\midrule
$c_2$ non-integral & $8$ & Proposition~\ref{prop:nonint}\\
$c_2$ integral, $a\ge\frac12$ & $8$ & Proposition~\ref{prop:largea}\\
$m\ge3$: $(4,1)$, $(3,1)$, $(3,2)$ & $8$ & Theorems~\ref{thm:41}, \ref{thm:31}, \ref{thm:32}\\
$m=2$, $K\ge3$: $(2,4)$, $(2,3)$ & $8$ & Theorems~\ref{thm:24}, \ref{thm:23}\\
$(2,1)$ & $8.486$ & Theorem~\ref{thm:21}\\
$(2,2)$ & $32/3$ & Theorem~\ref{thm:22}\\
$(1,1)$, $(1,2)$ & $32/3$ & Theorems~\ref{thm:11}, \ref{thm:12}\\
\bottomrule
\end{tabular}
\caption{Constants obtained for $\Gamma_{1,4}$ in each case of the canonical reduction ($c_2$ is the second shift coordinate, $a$ the minimum of the ternary form, $(m,K)$ as in Section~\ref{subsec:completion}). All other pairs $(m,K)$ are excluded by Proposition~\ref{prop:residual}. In the first four rows this paper gives independent proofs.}\label{tab:status}
\end{table}

\subsection*{Outline of the proof}
After the classical reductions (Section~\ref{sec:tools}) we may assume that $Q$ is in the canonical form
\[
Q=(x_1+a_2x_2+a_3x_3+a_4x_4+a_5x_5)x_2-\phi(x_3,x_4,x_5),
\]
where $\phi$ is a reduced positive definite ternary form with successive diagonal coefficients $a$, $A$, $t$ and $d:=(8|D|)^{1/5}$ satisfies $d^5=2aAt$. Completing the square in $x_3$, then in $x_4$, then in $x_5$ produces three successively simpler forms $F$, $G$, $H$ together with explicit target widths. The integers $m$ and $K$ record the quotient ranges governing the first two square-completion steps, and they index the residual cases. Section~\ref{sec:reduction} disposes of non-integral second shifts, of large minima $a\ge1/2$, and of all but nine pairs $(m,K)$:
\[
(4,1),\ (3,1),\ (3,2),\ (2,1),\ (2,2),\ (2,3),\ (2,4),\ (1,1),\ (1,2).
\]
Each of these is settled in Sections~\ref{sec:K1}--\ref{sec:m1}, independently of the case theorems of~\cite{BKR}. In particular Section~\ref{sec:reduction} gives independent proofs of the first two subcases of~\cite{BKR} mentioned above.

Three kinds of argument are combined.
\begin{enumerate}
\item \emph{Classical gap lemmas.} Dumir's square-completion lemma, Macbeath's lemma on gaps of $j+\beta z-\rho z^2$, and Jackson's interval lemma (Section~\ref{subsec:gaplemmas}). Applied with explicit width comparisons, they settle most of the parameter space and force the remaining obstructions onto finitely many \emph{resonances}: coefficients $q=H/N^2$ together with finitely many congruence classes of the accompanying linear coefficients and shifts (Lemma~\ref{lem:normal}).
\item \emph{Uniform residue gaps.} For a coefficient $q$ ranging over a rational interval, the largest circular gap of $\{\beta z-qz^2 \bmod 1: |z|\le M\}$, maximised over all $\beta$, can be computed exactly from a finite line arrangement (Proposition~\ref{prop:gapexact}). Such bounds settle whole ranges of $t$ at once. In particular a two-piece bound $\gap_4(t,\beta)\le\max(\frac12-t,4t)$ on $[\frac1{16},\frac16]$ (Lemma~\ref{lem:twopiece}) settles at once all branches with $L\ge4$ that occur, namely $(2,2)$ with $L=4$, $(1,2)$ with $L=4$, and $(2,3)$ with $L=4,5,6$ (Lemma~\ref{lem:highq}).
\item \emph{Exact finite certificates.} At a resonance the remaining problem is a family of affine inequalities in two or three real parameters; it is settled by a finite \emph{dyadic witness cover} (Definition~\ref{def:cover}). Families with many continuous parameters are settled by \emph{Farkas trees} over rational polytopes containing all admissible parameters (Definition~\ref{def:farkas}).
\end{enumerate}

\subsection*{The role of computation}
The proof combines analytic reductions with computer-verified universal covering statements. The analytic arguments, all given in the text, place every admissible parameter point in an explicitly specified domain; exact certificate checking then proves the required covering property on that domain. Each certificate is an explicit finite object: integer witnesses together with the boxes on which they satisfy given rational inequalities, or a tree of non-negative rational multipliers that combine the rows of an inequality system into a contradiction. The checkers use only exact rational arithmetic (Python \texttt{fractions}), so no floating-point number enters any acceptance decision. Floating-point linear programming is used only to \emph{discover} Farkas certificates, which are then accepted by exact replay. The verification data form a single public release with one entry point, described in Section~\ref{sec:verify}. It has two components:
\begin{itemize}
\item[(I)] $25$ Farkas trees ($23{,}753$ nodes) for the cases $(1,1)$, $(1,2)$ and $(2,1)$;
\item[(II)] $98$ dyadic witness covers, $3$ Farkas trees at the exact pinch points ($713$ nodes), exact residue-gap computations, two exact $t$-partitions and the ladder data.
\end{itemize}
The analytic reductions are not formalised in a proof assistant.

\subsection*{Organisation}
Section~\ref{sec:tools} collects the canonical form, the gap lemmas, square completion, residue gaps, resonant normal forms and the two certificate formats. Section~\ref{sec:reduction} reduces to the nine pairs $(m,K)$. Sections~\ref{sec:K1}--\ref{sec:m1} treat them, and Section~\ref{sec:end} completes the proof, gives the extremal form and describes the verification. Section~\ref{sec:concl} discusses which parts of the method extend beyond the present problem, and what remains open. Appendix~\ref{app:models} specifies the rational outer polytopes underlying the Farkas trees.
\section{Canonical form and tools}\label{sec:tools}

\subsection{The canonical form}\label{subsec:canonical}
Throughout, $Q$ has type $(1,4)$ and determinant $D\ne0$, and $d:=(8|D|)^{1/5}$. We must find, for each $c\in\R^5$, a point $X\in\Z^5+c$ with $0<Q(X)\le d$. The reduction to a canonical form uses the following classical inputs; a convenient account of their use for the present problem is~\cite[\S\S1--3]{BKR}.
\begin{enumerate}
\item[(I1)] Suppose that $Q$ is not a real multiple of a form with rational coefficients. Margulis's theorem on the Oppenheim conjecture~\cite{Margulis89} shows that $Q$ takes arbitrarily small nonzero values at integer points. Watson's theorem~\cite{Watson60} then shows that for every $\alpha\in\R$ and $\varepsilon>0$ the inequality $|Q(x+c)-\alpha|<\varepsilon$ has a solution $x\in\Z^5$, so $Q(\Z^5+c)$ is dense in $\R$. In particular it meets $(0,d)$, and there is nothing to prove.
\item[(I2)] If $Q$ is a multiple of a rational form, then it is isotropic over $\Q$ by Meyer's theorem (see~\cite[Ch.~IV, \S3.2]{Serre}), and Birch's reduction~\cite{Birch} transforms it, by an integral unimodular substitution and a scaling, into the canonical form below. Unimodular substitutions map the set $\Z^5+c$ to a set of the same kind, so the problem is unchanged.
\end{enumerate}
After these reductions,
\begin{align}
Q(X)&=(x_1+a_2x_2+a_3x_3+a_4x_4+a_5x_5)\,x_2-\phi(x_3,x_4,x_5),\label{eq:Q}\\
\phi(x,y,z)&=a(x+hy+kz)^2+A(y+\lambda z)^2+tz^2,\qquad d^5=8|D|=2aAt,\label{eq:phi}
\end{align}
where $a=\min_{0\ne w\in\Z^3}\phi(w)$, $|h|,|k|\le\tfrac12$, $0\le\lambda\le\tfrac12$ and $t+A\lambda^2\ge A$. Hermite's and Gauss's reduction inequalities (see~\cite{GL}) give
\begin{equation}\label{eq:red}
a^3\le d^5,\qquad A\ge\tfrac34a,\qquad t\ge\tfrac34A,\qquad 3aA^2\le 2d^5.
\end{equation}
We write $c=(c_1,\dots,c_5)$ for the shift. All shifts and all linear coefficients that are multiplied by an integer variable are significant only modulo~$1$. When $c_2\equiv0\pmod1$, an \emph{integral shear} $x_i\mapsto x_i+rx_2$ ($i\ge3$, $r\in\Z$) maps $\Z^5+c$ onto itself. It changes only the coefficients of $x_2x_i$ and $x_2^2$, so it may be used to normalise those coefficients. Reflections $x_i\mapsto -x_i$ ($i\ge3$) replace $(c_i,a_i)$ by $(-c_i,-a_i)$ and are used in the same way.

\subsection{Gap lemmas}\label{subsec:gaplemmas}
\begin{lemma}[Rounding]\label{lem:round}
Fix an admissible value $x_2=s\ne0$ and admissible values of $x_3,x_4,x_5$. As $x_1$ runs through $c_1+\Z$, the values of $Q$ form an arithmetic progression with difference $|s|$. Hence every interval $(u,u+|s|]$ contains a value. The same holds for the completed forms $F,G,H$ below.
\end{lemma}

If $c_2\equiv0$ we use $s=\pm1$. If $c_2\not\equiv0$ we let $s$ be the admissible value of least absolute value, so that $0<\eta:=|s|\le\frac12$. In that case only this one value of $s$ is used, never both signs.

\begin{lemma}[Dumir~\cite{Dumir68}; see {\cite[Lemma~3]{BKR}}]\label{lem:dumir}
Let $a>0$, $\gamma>a$, and let $m\ge1$ be the integer with $m<\gamma/a\le m+1$. Put $\delta=\gamma+\frac{m^2-1}{4}a$. Let $\theta\in\R$ and $F_0\in\R$ with $0<F_0\le\delta$. Then there is $x\in\theta+\Z$ with
\[
0<F_0+\tfrac a4-a x^2\le\gamma,
\]
and the upper inequality can be taken strict if $F_0<\delta$.
\end{lemma}

\begin{lemma}[Macbeath~\cite{Macbeath}; see {\cite[Lemma~5]{BKR}}]\label{lem:macbeath}
Let $\rho,\gamma>0$, $\beta\in\R$, $H\in\frac12\Z$ and $N\in\Z\setminus\{0\}$ with
\[
|H-N^2\rho|+\tfrac12<\gamma.
\]
Then for every $\nu\in\R$ there are $j,z\in\Z$ with $0<j+\nu+\beta z-\rho z^2<\gamma$, except possibly when
\[
\rho=H/N^2\quad\text{and}\quad \beta\in H/N+\Lambda,\qquad \Lambda=\tfrac1N\Z+2\rho\Z .
\]
\end{lemma}
We call $|H-N^2\rho|+\frac12$ the \emph{Macbeath width} of the pair $(H,N)$, and the exceptional value $\rho=H/N^2$ a \emph{resonance}. Only the strict form of the lemma is used.

\begin{lemma}[Jackson~\cite{Jackson}; see~\cite{BKR}]\label{lem:jackson}
Let $Q_0$ be a real quadratic form in $n_0$ variables with determinant $D_0\neq0$ that is a \emph{zero form}: $Q_0(w)=0$ for some nonzero $w\in\Z^{n_0}$ (equivalently, after clearing denominators, for some nonzero rational $w$). Let $u<v$ with $v-u>2|D_0|^{1/n_0}$. Then every shifted lattice $\Z^{n_0}+c'$ contains a point $x$ with $u<Q_0(x)<v$.
\end{lemma}
We apply it to the quadratic parts of $F$ (in $x_1,x_2,x_4,x_5$; $|D_0|=At/4$) and of $G$ (in $x_1,x_2,x_5$; $|D_0|=t/4$). Both vanish at the first unit vector, which is integral, and the constants $-a/4$, $-(a+A)/4$ are absorbed into $u$ and $v$.

\subsection{Square completion}\label{subsec:completion}
Completing the square in $x_3$, then $x_4$, then $x_5$ gives
\begin{align}
F&=(x_1+a_2'x_2+a_4'x_4+a_5'x_5)x_2-A(x_4+\lambda x_5)^2-tx_5^2-\tfrac a4,\label{eq:F}\\
G&=(x_1+a_2''x_2+a_5''x_5)x_2-tx_5^2-\tfrac{a+A}4,\label{eq:G}\\
H&=(x_1+a_2'''x_2)x_2-\tfrac{a+A+t}4,\label{eq:H}
\end{align}
with $a_2'=a_2+\frac{a_3^2}{4a}$, $a_4'=a_4-a_3h$, $a_5'=a_5-a_3k$, $a_2''=a_2'+\frac{(a_4')^2}{4A}$, $a_5''=a_5'-\lambda a_4'$ and $a_2'''=a_2''+\frac{(a_5'')^2}{4t}$. Indeed $Q=F+\frac a4-a(x_3+hx_4+kx_5-\frac{a_3x_2}{2a})^2$, and similarly for the other two steps. A shear $x_4\mapsto x_4+rx_2$ leaves $a_2''$ and $a_5''$ unchanged, and a shear $x_5\mapsto x_5+rx_2$ changes $a_5''$ by $-2tr$.

Assume $a<d$ and define integers $m,K,L\ge1$ and widths by
\begin{gather*}
m<\frac da\le m+1,\qquad \delta:=\delta_m=d+\frac{m^2-1}4a,\\
K<\frac{\delta}{A}\le K+1,\qquad \Delta:=\delta_{m,K}=\delta+\frac{K^2-1}4A,\\
L<\frac{\Delta}{t}\le L+1,\qquad \Delta_L:=\Delta+\frac{L^2-1}4t .
\end{gather*}
(The integers $K$ and $L$ are used only after it has been shown that $\delta>A$, respectively $\Delta>t$.) Three applications of Lemma~\ref{lem:dumir} give the following.

\begin{lemma}[Lifting]\label{lem:lift}
Fix any admissible $x_2=n\ne0$. If $0<H\le\Delta_L$ at some admissible $x_1$, then there is an admissible $x_5$ with $0<G\le\Delta$. If $0<G\le\Delta$, then there is an admissible $x_4$ with $0<F\le\delta$. If $0<F\le\delta$, then there is an admissible $x_3$ with $0<Q\le d$. Each implication remains true with strict upper inequalities.
\end{lemma}

By Lemma~\ref{lem:round} with $x_2=1$, the problem is solved if $\delta\ge1$, if $\Delta\ge1$ or if $\Delta_L\ge1$ (when $c_2\equiv0$), and with $1$ replaced by $\eta$ when $c_2\not\equiv0$. If $t\ge\Delta/(L+1)$, then
\begin{equation}\label{eq:DeltaL}
\Delta_L\ge\tfrac{L+3}{4}\Delta .
\end{equation}

\subsection{Uniform residue gaps}\label{subsec:gaps}
\begin{definition}\label{def:gap}
For $q,\beta\in\R$ and $M\ge1$, let $\gap_M(q,\beta)$ be the largest gap between cyclically consecutive points of $\{\beta z-qz^2\bmod 1:\ |z|\le M\}\subset\R/\Z$. For a closed interval $I$ put $\gap_M(I)=\sup_{q\in I,\beta\in\R}\gap_M(q,\beta)$.
\end{definition}

\begin{lemma}\label{lem:gapuse}
If $\gap_M(q,\beta)<W$, then for every $\nu\in\R$ there are $j\in\Z$ and $|z|\le M$ with $0<j+\nu+\beta z-qz^2<W$.
\end{lemma}
\begin{proof}
The open arc $(-\nu,-\nu+W)$ of $\R/\Z$ has length larger than every gap, so it contains one of the points.
\end{proof}

\begin{proposition}\label{prop:gapexact}
For rational $q_-\le q_+$, the number $\gap_M([q_-,q_+])$ is attained, and it equals the maximum of $\gap_M$ over the finitely many points of $[0,\frac12]\times[q_-,q_+]$ at which two of the following lines meet: the boundary lines of this rectangle, and the collision lines $(i-j)\beta-(i^2-j^2)q=k$ with $-M\le i<j\le M$ and $k\in\Z$.
\end{proposition}
\begin{proof}
Replacing $z$ by $-z$ changes $\beta$ to $-\beta$, and $\beta$ matters modulo $1$, so we may take $0\le\beta\le\frac12$. Two residues coincide exactly on a collision line, and only finitely many collision lines meet the rectangle. On each closed cell of the resulting arrangement the cyclic order of the residues is fixed, so every consecutive gap is an affine function there. Hence $\gap_M$ is a maximum of finitely many affine functions on each cell. It is therefore convex on the cell and attains its maximum at a vertex. The function $\gap_M$ is continuous, so the cell maxima give the global maximum.
\end{proof}

With $x_2=1$ and $x_5=z+c_5$, the values of $G$ are $j+\nu+\beta z-tz^2$ with $\beta\equiv a_5''-2tc_5$ and $\nu$ depending on the other data. Thus $\gap_M(\{t\})<\Delta$ gives $0<G<\Delta$. Likewise, with $x_2=1$ and $x_5$ fixed, the values of $F$ are $j+\nu+\beta y-Ay^2$, and $\gap_M(\{A\})<\delta$ gives $0<F<\delta$.

\begin{lemma}[Two-piece gap bound]\label{lem:twopiece}
For every real $\beta$,
\[
\gap_4(t,\beta)\le\begin{cases}\frac12-t,&\frac1{16}\le t\le\frac1{10},\\[2pt] 4t,&\frac1{10}\le t\le\frac16.\end{cases}
\]
\end{lemma}
\begin{proof}
Use the arrangement of Proposition~\ref{prop:gapexact} on the rectangles $[0,\frac12]\times[\frac1{16},\frac1{10}]$ and $[0,\frac12]\times[\frac1{10},\frac16]$. On each cell every consecutive gap is affine, and so is the bound, so $\gap_4$ minus the bound is convex on each cell. It therefore suffices to check the vertices. Exact rational arithmetic verifies the bound at all $160$ and $229$ vertices, respectively.
\end{proof}

\begin{lemma}[High-quotient closure]\label{lem:highq}
Suppose $4\le L\le6$ and
\[
Lt<\Delta\le(L+1)t,\qquad \Delta+\tfrac{L^2-1}4t\le1,\qquad \Delta+t>\tfrac12 .
\]
Then for all $\beta,\nu\in\R$ there are $j,z\in\Z$ with $|z|\le4$ and $0<j+\nu+\beta z-tz^2<\Delta$. In particular, if these conditions hold for the parameters of $G$, then $0<G<\Delta$ with $x_2=1$.
\end{lemma}
\begin{proof}
From $(L+2)t\ge\Delta+t>\frac12$ we get $t>\frac1{2(L+2)}\ge\frac1{16}$. From $Lt+\frac{L^2-1}4t<\Delta+\frac{L^2-1}4t\le1$ we get $t<\frac4{L^2+4L-1}\le\frac4{31}<\frac16$. By Lemma~\ref{lem:twopiece}, $\gap_4(t,\beta)\le\max(\frac12-t,4t)$. Here $\frac12-t<\Delta$ by the third hypothesis and $4t\le Lt<\Delta$, so Lemma~\ref{lem:gapuse} applies.
\end{proof}

\subsection{Resonant normal forms}\label{subsec:normal}
\begin{lemma}\label{lem:normal}
Let $c_2\equiv0$. Consider $G$ with $x_2=\pm1$, the variable $x_5=z+c$ ($c=c_5$), the coefficient $q=t$ and $e=a_5''$. Alternatively consider $F$ with $x_2=\pm1$ and $x_5\in\{c_5,c_5+1\}$, the variable $x_4=y+c$ ($c=c_4$), $q=A$ and $e=a_4'$. Let $(H,N)$ be a row of Table~\ref{tab:normal} whose Macbeath width is smaller than the target; in the case of $G$ any row may be used, in the case of $F$ only the first three rows. Suppose Lemma~\ref{lem:macbeath} gives no value in the target interval for any of these choices. Then $q=H/N^2$, and in the case of $F$ also $\lambda=0$. After an integral shear of the variable by $x_2$ and possibly a reflection of it, $(e,c)$ is congruent modulo~$1$ to one of the normal forms in the table.
\end{lemma}

\begin{table}[ht]
\centering
\begin{tabular}{cccl}
\toprule
$(H,N)$ & $q=H/N^2$ & $\Lambda$ & normal forms $(e,c)$\\
\midrule
$(1/2,1)$ & $1/2$ & $\Z$ & $(0,1/2),\ (1/2,0)$\\
$(1,2)$ & $1/4$ & $\frac12\Z$ & $(0,0),\ (1/4,1/2)$\\
$(3/2,3)$ & $1/6$ & $\frac13\Z$ & $(1/6,0),\ (0,1/2)$\\
$(5,5)$ & $1/5$ & $\frac15\Z$ & $(0,0),\ (1/10,1/4),\ (1/5,1/2)$\\
\bottomrule
\end{tabular}
\caption{Resonances and normal forms (Lemma~\ref{lem:normal}).}\label{tab:normal}
\end{table}

\begin{proof}
In each row $2q\in\frac1N\Z$, so $\Lambda=\frac1N\Z$. With $x_2=\pm1$ the linear coefficient of $z$ is $\beta_\pm\equiv\pm e-2qc$, and failure for both signs means $\beta_\pm\in H/N+\Lambda$. Adding gives $4qc\in -2H/N+\Lambda$, and subtracting gives $2e\in\Lambda$. A shear changes $e$ by multiples of $2q$, a reflection maps $(e,c)$ to $(-e,-c)$, and $e$ matters modulo~$1$.
\begin{itemize}
\item $(1/2,1)$: $2c\in\Z$, and $e\equiv\frac12+c\pmod 1$.
\item $(1,2)$: $c\in\frac12\Z$ and $e\equiv\frac c2\pmod{\frac12}$. Shears by $2q=\frac12$ give $e=\frac c2$.
\item $(3/2,3)$: $2c\in\Z$ and $e\equiv\frac12+\frac c3\pmod{\frac13}$. Shears by $\frac13$ give the two listed forms.
\item $(5,5)$: $4c\in\Z$ and $e\equiv\frac{2c}5\pmod{\frac15}$. Shears move $e$ by $\frac{2r}5$, which runs through all of $\frac15\Z/\Z$, so we may take $e=\frac{2c}5$. Finally the reflection identifies $c=\frac34$ with $c=\frac14$.
\end{itemize}
For $F$, replacing $x_5$ by $x_5+1$ changes the linear coefficient of $y$ by $-2A\lambda$. Failure for both values forces $2A\lambda=2H\lambda/N^2\in\frac1N\Z$, that is $2H\lambda/N\in\Z$. In the three rows allowed for $F$ ($H/N=\frac12$) this gives $\lambda\in\Z$, so $\lambda=0$ as $0\le\lambda\le\frac12$. (For the row $(5,5)$ it would give only $2\lambda\in\Z$, which is why that row is excluded for $F$.)
\end{proof}

\subsection{Exact certificates}\label{subsec:certs}
\begin{definition}[Fixed family and dyadic cover]\label{def:cover}
A \emph{fixed family} consists of the following data:
\begin{itemize}
\item for each integer $n\ne0$, a set $W_n$ of \emph{admissible labels} $w$ (in all our families, integer vectors $w$, all of which are admissible), with rational constants $K(n,w)$;
\item a rational interval $[\mu_-,\mu_+]$;
\item a finite list of \emph{targets} $(c_1,c_0)\in\Q^2$.
\end{itemize}
Its witness values are
\[
f(n,w,j;\alpha,\kappa,\mu)=nj+n^2\alpha+n\kappa+K(n,w)-\mu\qquad(n\ne0,\ w\in W_n,\ j\in\Z).
\]
The family itself permits all admissible witnesses $(n,w,j)$. A \emph{finite certificate} records a bounded subset of them, with $|n|\le n_{\max}$ and finitely many labels; a \emph{ladder} (Section~\ref{subsec:pinch}) records an explicit infinite subfamily.

A \emph{dyadic cover} is a finite certificate organised as a finite rooted tree of boxes. The root is $[0,1]^2\times[\mu_-,\mu_+]$ in $(\alpha,\kappa,\mu)$. Each internal node is split into its four quadrants in $(\alpha,\kappa)$ or into two halves in $\mu$. Each leaf $B$ either carries a label $(n,w,j)$ such that $\min_Bf>0$ and $\max_B(f-c_1\mu)<c_0$ for some target, or is a \emph{pinch leaf}: its $(\alpha,\kappa)$-projection lies, modulo $\Z^2$, in a closed $\ell^1$-ball $\bar B(p^*,r)$ about a point $p^*$ for which Lemma~\ref{lem:ladder} below holds, for every $\mu$ in the $\mu$-range of the leaf, with some radius larger than $r$.
\end{definition}

\begin{lemma}\label{lem:cover}
Suppose a fixed family has a dyadic cover, each of whose targets satisfies $c_0+c_1\mu\le W$ for the relevant width $W$ on the whole domain. Let $p=(\alpha,\kappa,\mu)\in\R^2\times[\mu_-,\mu_+]$. Then there are $n\ne0$, $w$, $j$ with $0<f(n,w,j;p)<W$, unless $(\alpha,\kappa)\equiv p^*\pmod{\Z^2}$ for the centre $p^*$ of a pinch ball.
\end{lemma}
\begin{proof}
Replacing $\alpha$ by $\alpha+1$ or $\kappa$ by $\kappa+1$ changes $f$ by $n^2$ or $n$, which can be absorbed into $j$. So we may assume $(\alpha,\kappa)\in[0,1]^2$, and $p$ lies in some leaf $B$. If $B$ is a pinch leaf, Lemma~\ref{lem:ladder} applies. Otherwise $f$ is affine with coefficients $n^2>0$, $n$ and $-1$ in $\alpha$, $\kappa$ and $\mu$. Its extrema over $B$ are therefore at explicit corners, and the leaf inequalities hold at $p$.
\end{proof}

\begin{proposition}[Completeness of dyadic covers]\label{prop:complete}
Let $D$ be the root box of a fixed family with finitely many targets. Let $U=\bigcup_{i=1}^sB(p_i,r_i)$ be a finite union of open $\ell^1$-balls in the $(\alpha,\kappa)$-torus, where for each $i$ Lemma~\ref{lem:ladder} holds at $p_i$, for every $\mu\in[\mu_-,\mu_+]$, with some radius $R_i>r_i$. Suppose every point of $D$ whose $(\alpha,\kappa)$ lies outside $U$ admits some admissible witness $(n,w,j)$ of the family and some target with $0<f<c_0+c_1\mu$. Then there is a dyadic cover whose witnesses are finitely many, and each of whose pinch leaves projects into one of the closed balls $\bar B(p_i,\rho_i)$, for some fixed $\rho_i$ with $r_i<\rho_i<R_i$.
\end{proposition}
\begin{proof}
Fix $\rho_i\in(r_i,R_i)$. For a witness $\omega=(n,w,j)$ and a target $\tau=(c_1,c_0)$, the set $O_{\omega,\tau}$ of points of $D$ with $f>0$ and $f-c_1\mu<c_0$ is open in $D$. The set $K$ of points of $D$ whose $(\alpha,\kappa)$ lies outside $U$ is compact, so finitely many sets $O_{\omega,\tau}$ cover it. Let $V_i$ be the set of points of $D$ whose $(\alpha,\kappa)$ lies in $B(p_i,\rho_i)$ modulo $\Z^2$; it is open in $D$. The finitely many sets $O_{\omega,\tau}$ together with $V_1,\dots,V_s$ form an open cover of the compact metric space $D$, since every point outside $K$ lies in some $V_i$. Let $\lambda>0$ be a Lebesgue number for this cover, and subdivide $D$ dyadically until every box has diameter less than $\lambda$. Each closed box $B$ then lies in one member of the cover. If $B\subset O_{\omega,\tau}$, then by compactness $\min_Bf>0$ and $\max_B(f-c_1\mu)<c_0$, so $B$ is a witness leaf. If $B\subset V_i$, its projection lies in $\bar B(p_i,\rho_i)$ and $R_i>\rho_i$, so $B$ is a pinch leaf.
\end{proof}

Thus compactness supplies the finite part of a certificate, while each ladder ball must be certified separately; the ladders are the only infinite ingredient.

\subsection{Pinch points}\label{subsec:pinch}
For a fixed family with $w=z\in\Z$ and $K(n,z)=ne(z+c)-q(z+c)^2$, and $p=(\alpha,\kappa)$, let
\[
V_n(p)=\{nj+n^2\alpha+n\kappa+K(n,z):\ j,z\in\Z\},\qquad V(p)=\textstyle\bigcup_{n\ne0}V_n(p).
\]
Then $V_n(p^*+\varepsilon)=V_n(p^*)+n^2\varepsilon_1+n\varepsilon_2$. We call $(v,n)$ with $n\ge1$ a \emph{rung} at $p^*$ if $v\in V_n(p^*)\cap V_{-n}(p^*)$. A \emph{ladder} is a sequence of rungs $(v_k,n_k)_{k\ge1}$ with $1=n_1<n_2<\cdots$, $n_k\to\infty$ and finitely many distinct values $v_k$. Put $\rho_k=n_{k+1}/n_k$.

\begin{lemma}[Ladder lemma]\label{lem:ladder}
Let $\mu,\Delta>0$. Suppose that $p^*$ has an ascending ladder $(v_k^-,n_k^-)$ and a descending ladder $(v_k^+,n_k^+)$ with
\[
\text{(A)}\quad \rho_k^2\,(\mu-v_k^-)<\mu+\Delta-v_{k+1}^-,\qquad \text{(D)}\quad \rho_k^2\,(v_k^+-\mu-\Delta)<v_{k+1}^+-\mu\qquad(k\ge1).
\]
Put $r=\min(v_1^+-\mu,\ \mu+\Delta-v_1^-)$. Then for every $\varepsilon\in\R^2$ with $0<|\varepsilon_1|+|\varepsilon_2|<r$ the set $V(p^*+\varepsilon)$ meets $(\mu,\mu+\Delta)$.
\end{lemma}
\begin{proof}
Suppose first $\varepsilon_1\ge0$. Let $s=\operatorname{sign}\varepsilon_2$ (or $s=1$ if $\varepsilon_2=0$) and $g_k=n_k^2\varepsilon_1+n_k|\varepsilon_2|$, where $n_k=n_k^-$. Then $u_k=v_k^-+g_k\in V_{sn_k}(p^*+\varepsilon)$, $g_k$ increases to infinity, and $g_{k+1}\le\rho_k^2g_k$. If $u_1>\mu$, then $u_1<v_1^-+r\le\mu+\Delta$. Otherwise let $k$ be the largest index with $u_k\le\mu$. Then $u_{k+1}>\mu$ and $g_{k+1}\le\rho_k^2(\mu-v_k^-)<\mu+\Delta-v_{k+1}^-$ by (A), so $u_{k+1}<\mu+\Delta$.

If $\varepsilon_1<0$, let $s=-\operatorname{sign}\varepsilon_2$ (or $1$), $h_k=n_k^2|\varepsilon_1|+n_k|\varepsilon_2|$ with $n_k=n_k^+$, and $u_k=v_k^+-h_k\in V_{sn_k}(p^*+\varepsilon)$. If $u_1<\mu+\Delta$, then $u_1>v_1^+-r\ge\mu$. Otherwise let $k$ be the largest index with $u_k\ge\mu+\Delta$. Then $u_{k+1}<\mu+\Delta$ and $h_{k+1}\le\rho_k^2(v_k^+-\mu-\Delta)<v_{k+1}^+-\mu$ by (D), so $u_{k+1}>\mu$.
\end{proof}

Both (A) and (D) improve as $\Delta$ increases. So it suffices to verify them for a lower bound for $\Delta$ that is affine in $\mu$, at the two endpoints of the $\mu$-range. When a value lies in the window at $p^*$ itself, the point $p^*$ is not a failure point either. Otherwise $p^*$ is an isolated failure point: a \emph{pinch point}.

\begin{definition}[Farkas tree]\label{def:farkas}
Let $P=\{\theta:B_i\theta\le b_i\}$ be a rational polytope and let $T(\theta)$ be an affine target. A \emph{Farkas tree} for $(P,T)$ is a finite binary tree whose internal nodes carry integer witnesses $\omega$ with affine values $f_\omega(\theta)$. The two children of a node add the rows $f_\omega\le0$ and $f_\omega\ge T$ respectively; for the non-strict goal $0<f\le T$ the second row is instead $f_\omega>T$. Each leaf carries rational multipliers $w_i\ge0$ on the rows of its path with $\sum w_iB_i=0$, and either $\sum w_ib_i<0$, or $\sum w_ib_i=0$ with positive weight on a strict row.
\end{definition}

\begin{lemma}\label{lem:farkas}
If $(P,T)$ has a Farkas tree, then for every $\theta\in P$ some witness on the path through $\theta$ satisfies $0<f_\omega(\theta)<T(\theta)$ (respectively $\le T(\theta)$).
\end{lemma}
\begin{proof}
Follow the branch whose added rows hold at $\theta$. Such a branch cannot end at a leaf, since the leaf certifies that its rows are inconsistent. So at some node neither failure row holds, and that node's witness succeeds at $\theta$.
\end{proof}

The polytopes are \emph{outer} relaxations: Appendix~\ref{app:models} lists inequalities that every admissible parameter point satisfies, written in moment coordinates in which every witness is affine. Equality with the admissible set is never needed.
\section{Reduction to nine residual cases}\label{sec:reduction}

\subsection{Quotient pruning}
The following criterion holds in any dimension, so we state it for canonical forms of type $(1,n-1)$:
\[
Q=(x_1+\ell_2x_2+\cdots+\ell_nx_n)x_2-\Phi(x_3,\dots,x_n),\qquad \Phi>0,\qquad d^n=C|D|.
\]
Suppose successive square completions remove the coefficients $a_1,\dots,a_k$ of $\Phi$, with widths
\[
W_0=d,\qquad m_i<\frac{W_{i-1}}{a_i}\le m_i+1,\qquad W_i=W_{i-1}+\frac{m_i^2-1}{4}a_i\qquad(m_i\ge1).
\]
For $n=5$ and $C=8$ these are $a_1=a$, $a_2=A$, $a_3=t$ and $W_1=\delta$, $W_2=\Delta$, $W_3=\Delta_L$ of Section~\ref{subsec:completion}.

\begin{lemma}[Uniform quotient pruning]\label{lem:prune}
Let $1\le k\le n-3$. If
\begin{equation}\label{eq:prune}
\prod_{i=1}^{k}\frac{(m_i+3)^{\,n-i}}{m_i+1}>\frac{2^{\,n-k}}{C}\,4^{\,kn-k(k+1)/2},
\end{equation}
then every shift class contains $X$ with $0<Q(X)<d$.
\end{lemma}
\begin{proof}
Since $a_i\ge W_{i-1}/(m_i+1)$, we have $W_i\ge W_{i-1}(1+\frac{m_i-1}4)=\frac{m_i+3}4W_{i-1}$. Hence $W_{i}\ge d\prod_{l\le i}\frac{m_l+3}4$, and
\[
W_k^{\,n-k}\prod_{i=1}^ka_i\ \ge\ W_k^{\,n-k}\prod_{i=1}^k\frac{W_{i-1}}{m_i+1}\ \ge\ d^n\prod_{i=1}^k\frac{((m_i+3)/4)^{\,n-i}}{m_i+1}>2^{\,n-k}|D|
\]
by~\eqref{eq:prune}. After the $k$ completions, the homogeneous part of the remaining form is isotropic in $n-k$ variables with determinant of absolute value $|D|/(a_1\cdots a_k)$. So $W_k>2(|D|/(a_1\cdots a_k))^{1/(n-k)}$, and Lemma~\ref{lem:jackson} gives a value of the remaining form in $(0,W_k)$. Its $x_2$-coordinate is nonzero, since for $x_2=0$ the form is negative. The $k$ completions are then undone by Lemma~\ref{lem:dumir}, as in Lemma~\ref{lem:lift}.
\end{proof}

For every $p\ge3$ the function $(m+3)^p/(m+1)$ increases for $m\ge1$, since its logarithmic derivative is $\frac{(p-1)m+p-3}{(m+1)(m+3)}\ge0$. So a single check at the smallest value of each quotient eliminates all larger values. For $(n,C)=(5,8)$, condition~\eqref{eq:prune} reads $(m+3)^4>512(m+1)$ when $k=1$ and $(m+3)^4(K+3)^3>4^7(m+1)(K+1)$ when $k=2$.

\begin{corollary}\label{lem:small}
If $a<d/5$, then every shift class contains $X$ with $0<Q(X)<d$.
\end{corollary}
\begin{proof}
Here $m\ge5$, and $(5+3)^4=4096>512\cdot6$.
\end{proof}

\subsection{Non-integral second shift}
\begin{proposition}\label{prop:nonint}
If $c_2\not\equiv0\pmod1$, then every shift class contains $X$ with $0<Q(X)\le d$.
\end{proposition}
\begin{proof}
Let $s$ be the admissible value of $x_2$ of least absolute value, $\eta=|s|\le\frac12$. If $d\ge\eta$, Lemma~\ref{lem:round} applies. So let $d<\eta\le\frac12$. Then $a\le d^{5/3}<d$, and by Corollary~\ref{lem:small} we may assume $1\le m\le4$.

\emph{Case $m=1$, i.e.\ $d/2\le a<d$.} From $a^3\le d^5$ we get $d\ge8^{-1/2}$, hence $a+d\ge\frac32d>\frac12\ge\eta$. Divide $Q$ by $\eta$ and apply Lemma~\ref{lem:macbeath} in $x_3$ with $(H,N)=(\frac12,1)$. The width condition $|\frac12-\frac a\eta|+\frac12<\frac d\eta$ reads $\max(a,\eta-a)<d$, which holds. Failure forces $a=\eta/2$. Changing the integer parts of $x_4,x_5$ changes the linear coefficient of $x_3$ by $2ah/\eta=h$ and by $k$, so failure also forces $h=k=0$. Then $A\ge a=\eta/2$. Moreover $A\le\frac2{\sqrt3}d^2<d$ by~\eqref{eq:red}, so the same argument in $x_4$ forces $A=\eta/2$ and $\lambda=0$. Then $t\ge A=\eta/2$ and $d^5=2aAt\ge\eta^3/4>d^3/4$, contradicting $d<\frac12$.

\emph{Case $2\le m\le4$.} By Lemma~\ref{lem:lift} with $x_2=s$ and Lemma~\ref{lem:round}, we may assume $\delta<\eta\le\frac12$, so $\frac{m+3}4d\le\delta<\frac12$. From $A^2\le 2d^5/(3a)\le\frac23(m+1)d^4$ we get
\[
\frac{\delta}{A}\ge\frac{m+3}{4d}\sqrt{\frac{3}{2(m+1)}}>\frac{(m+3)^2}{8}\sqrt{\frac{3}{2(m+1)}} ,
\]
which exceeds $3$ for $m=4$ and $2$ for $m=2,3$. Hence $K\ge3$ for $m=4$ and $K\ge2$ for $m=2,3$. We have $a\ge\frac{d}{m+1}$, $A\ge\frac{\delta}{K+1}\ge\frac{(m+3)d}{4(K+1)}$ and $\Delta\ge\frac{K+3}4\delta\ge\frac{(m+3)(K+3)d}{16}$.
\begin{enumerate}
\item Lemma~\ref{lem:prune} with $k=2$ solves the problem when
\begin{equation}\label{eq:jackson7}
(m+3)^4(K+3)^3>4^7(m+1)(K+1).
\end{equation}
This holds for $(m,K)=(4,2),(3,3),(2,5)$, hence for $(4,K\ge2)$, $(3,K\ge3)$ and $(2,K\ge5)$.
\item $(m,K)=(3,2)$: $A\ge\delta/3\ge d/2$ and $A^2\le\frac83d^4$ give $d^2\ge\frac3{32}$. Hence $\Delta^2\ge(\frac{15}8)^2d^2\ge\frac{675}{2048}>\frac14$, so $\Delta>\frac12\ge\eta$ and rounding solves $G$.
\item $m=2$, $K\in\{2,3,4\}$: if $\Delta<\eta$ then $d<\frac{8}{5(K+3)}$. With $a\ge\frac d3$ and $A\ge\frac{5d}{4(K+1)}$,
\[
\frac\Delta t=\frac{2aA\Delta}{d^5}\ge\frac{25(K+3)}{96(K+1)d^2}>\frac{625(K+3)^3}{6144(K+1)}\ge\frac{78125}{18432}>4 .
\]
So $L\ge4$, and by~\eqref{eq:DeltaL} $\Delta_L\ge\frac74\Delta\ge\frac{175}{64}d\ge\frac{175}{64\sqrt{27}}>\frac12\ge\eta$, using $d\ge27^{-1/2}$, which follows from $a\ge d/3$ and $a^3\le d^5$. Rounding solves $H$. \qedhere
\end{enumerate}
\end{proof}

\subsection{Integral second shift and large minimum}
From now on $c_2\equiv0$, so both $x_2=1$ and $x_2=-1$ are admissible, and integral shears are available.

\begin{proposition}\label{prop:largea}
If $c_2\equiv0$ and $a\ge\frac12$, then every shift class contains $X$ with $0<Q(X)\le d$. If $d<1$, then even $0<Q(X)<d$.
\end{proposition}
\begin{proof}
If $d\ge1$, rounding gives $0<Q\le1\le d$. So let $d<1$; then $a\le d^{5/3}<d$. If $a>\frac12$, Lemma~\ref{lem:macbeath} in $x_3$ with $(H,N)=(\frac12,1)$ has width $a<d$ and no resonance. So let $a=\frac12$. Then $d^5\ge a^3=\frac18$. Failure of the same argument for $x_2=\pm1$ and all integer parts of $x_4,x_5$ forces $h=k=0$ and $(a_3,c_3)\equiv(0,\frac12)$ or $(\frac12,0)$, by the argument of Lemma~\ref{lem:normal} applied to $Q$ in $x_3$. Since $h=0$, also $A\ge a=\frac12$.

\emph{If $A<d$:} the same argument in $x_4$ forces $A=\frac12$, $\lambda=0$ and $(a_4,c_4)\in\{(0,\frac12),(\frac12,0)\}$. Then $t=2d^5\in[\frac12,2)$ and $d\ge4^{-1/5}>\frac34$. Let $H$ be a half-integer nearest to $t$. The width with $(H,1)$ is at most $\frac34<d$, so failure forces $t=r/2$ with $r\in\{1,2,3\}$ and $\pm a_5-rc_5\equiv\frac r2\pmod1$.

\emph{If $A\ge d$:} then $d^2\le A^2\le\frac43d^5$, so $d\ge(\frac34)^{1/3}>\frac34$ and $A<\frac54$. The pair $(1,1)$ in $x_4$ has width $\frac32-A\le\frac32-d<d$ if $A\le1$, and $A-\frac12<\frac34<d$ if $A\ge1$. Failure forces $A=1$ and $2\lambda\in\Z$. If $\lambda=0$ then $t\ge1>d^5=t$, which is impossible, so $\lambda=\frac12$. Then $C:=t+\frac14\in[1,\frac54)$ is the coefficient of $x_5^2$. The pair $(1,1)$ in $x_5$ has width $C-\frac12<d$, so failure forces $C=1$, i.e.\ $t=\frac34$.

The remaining negative parts are therefore
\[
\phi_r=\tfrac12(x^2+y^2+rz^2)\quad(r=1,2,3,\ d^5=r/4),\qquad \phi_H=\tfrac12x^2+y^2+yz+z^2\quad(d^5=\tfrac34).
\]
Write $b=(a_3,a_4,a_5)$ and $c=(c_3,c_4,c_5)$, indexed by the three negative coordinates. Using both signs of $x_2$ in the Macbeath failures, the admissible classes are as follows. For $\phi_r$: $(b_1,c_1),(b_2,c_2)\in\{(0,\frac12),(\frac12,0)\}$, $b_3\in\{0,\frac12\}$ and $c_3\equiv(b_3-\frac r2+j)/r$ with $0\le j<r$, giving $8r$ classes. For $\phi_H$: $(b_1,c_1)$ as before, $b_2,b_3\in\{0,\frac12\}$, $2c_2+c_3\equiv b_2$ and $c_2+2c_3\equiv b_3$, giving three solutions for each $(b_2,b_3)$ and $24$ classes. With $x_2=n$, $(x_3,x_4,x_5)=w+c$, $x_1=j+c_1^{\rm orig}$, $\alpha=a_2$ and $\kappa=c_1^{\rm orig}+b\cdot c$,
\[
Q=nj+n^2\alpha+n\kappa+n\,b\cdot w-\phi(w+c).
\]
This is a fixed family with $K(n,w)=n\,b\cdot w-\phi(w+c)$, $w\in\{-1,0,1\}^3$, $|n|\le4$, $\mu=0$ and the target $\Delta_\phi$ of Table~\ref{tab:largea}. Each of the $72$ families has a dyadic cover, so by Lemma~\ref{lem:cover} $0<Q<\Delta_\phi$. Finally $\Delta_\phi<d$, since $(\frac{151}{200})^5<\frac14$, $(\frac{87}{100})^5<\frac12$ and $(\frac{15}{16})^5<\frac34$.
\end{proof}

\begin{remark}
The strict inequality can fail when $d=1$. For $\phi=x^2+y^2+z^2+xy+xz+yz$ one has $a=1$, $h=k=\frac12$, $A=\frac34$, $\lambda=\frac13$, $t=\frac23$ and $d^5=2aAt=1$. The form $Q=x_1x_2-\phi$ is integral, so it has no value in $(0,1)$ on $\Z^5$, while $Q(1,1,0,0,0)=1$.
\end{remark}

\begin{table}[ht]
\centering
\begin{tabular}{lcr}
\toprule
Block & $\Delta_\phi$ & Classes\\
\midrule
$\phi_1$ & $151/200$ & 8\\
$\phi_2$ & $87/100$ & 16\\
$\phi_3$ & $15/16$ & 24\\
$\phi_H$ & $15/16$ & 24\\
\bottomrule
\end{tabular}
\caption{Dyadic covers for the forced classes with $a=\frac12$ (Proposition~\ref{prop:largea}).}\label{tab:largea}
\end{table}

\subsection{The residual cases}
\begin{proposition}\label{prop:residual}
Assume $c_2\equiv0$. The theorem holds unless
\begin{equation}\label{eq:residual}
d<1,\qquad a<\tfrac12,\qquad a+d\le1,
\end{equation}
and $(m,K)$ is one of $(4,1)$, $(3,1)$, $(3,2)$, $(2,1)$, $(2,2)$, $(2,3)$, $(2,4)$, $(1,1)$, $(1,2)$. In all these cases $\delta<1$ may also be assumed.
\end{proposition}
\begin{proof}
Rounding handles $d\ge1$ (with $0<Q\le d$) and Proposition~\ref{prop:largea} handles $a\ge\frac12$. If $a<\frac12$ and $a+d>1$, the pair $(\frac12,1)$ in $x_3$ has width $1-a<d$ and no resonance. Under~\eqref{eq:residual}, $a\le d^{5/3}<d$, and Corollary~\ref{lem:small} leaves $1\le m\le4$. Rounding handles $\delta\ge1$; otherwise
\[
a\ge\frac d{m+1},\qquad \delta\ge\frac{m+3}4d,\qquad d\le B_m:=\min\Bigl(\frac{m+1}{m+2},\frac4{m+3}\Bigr),
\]
and $\delta/A\ge\frac{m+3}{4d}\sqrt{3/(2(m+1))}>1$ for $1\le m\le4$, by inserting $d\le B_m$. Thus $K\ge1$. Lemma~\ref{lem:prune} settles every pair satisfying~\eqref{eq:jackson7}, which leaves $K=1$ for $m=4$, $K\le2$ for $m=3$, and $K\le4$ for $m=2$. For $m=1$ we have $\delta=d$ and $d/A\le\frac{4d}{3a}\le\frac83$, so $K\le2$.
\end{proof}

Exponent three on $K+3$ in~\eqref{eq:jackson7} is correct (it is $n-k$ in~\eqref{eq:prune}). The fourth power printed in~\cite[(4.3), (6.4)]{BKR} is a misprint, as comparison with~\cite[(3.17)]{BKR} shows. The conventions $m<d/a\le m+1$ and $K<\delta/A\le K+1$ assign every boundary point to exactly one case; certified domains may contain additional boundary points. Table~\ref{tab:cases} lists where each residual pair is treated.

\begin{table}[ht]
\centering
\begin{tabular}{lll}
\toprule
$(m,K)$ & Section & Method\\
\midrule
$(4,1)$ & \ref{sec:K1} & analytic\\
$(3,1)$ & \ref{sec:K1} & interval tests, residue gaps, 1 dyadic cover\\
$(3,2),(2,4),(2,3)$ & \ref{sec:third} & Macbeath, uniform gaps, high-quotient lemma, 11 covers\\
$(2,2)$ & \ref{sec:22} & high-quotient lemma, $t$-partitions, 10 covers, 3 corner trees\\
$(2,1)$ & \ref{sec:21} & ternary and diagonal Farkas trees (component I)\\
$(1,1),(1,2)$ & \ref{sec:m1} & Farkas trees (component I), high-quotient lemma, 4 covers\\
\bottomrule
\end{tabular}
\caption{The nine residual pairs.}\label{tab:cases}
\end{table}
\section{\texorpdfstring{The cases $(4,1)$ and $(3,1)$}{The cases (4,1) and (3,1)}}\label{sec:K1}

Throughout this section $c_2\equiv0$, \eqref{eq:residual} holds, $m\in\{3,4\}$, $K=1$ and $\delta<1$. Thus
\begin{equation}\label{eq:K1dom}
\frac d{m+1}\le a<\frac dm,\qquad \frac\delta2\le A<\delta,\qquad \delta\ge\frac{m+3}4d,\qquad d<\frac4{m+3},
\end{equation}
and $3aA^2\le2d^5$ gives $d^2\ge\frac{3(m+3)^2}{128(m+1)}$, i.e.\ $d^2\ge\frac{27}{128}$ for $m=3$ and $d^2\ge\frac{147}{640}$ for $m=4$. In both cases $\delta>\frac23$.

\begin{lemma}[Diagonal reduction]\label{lem:diag}
If $(m,K)\in\{(3,1),(4,1)\}$ and no $X\equiv c$ has $0<Q(X)\le d$, then
\[
A=t=\tfrac12,\qquad \lambda=h=k=0,\qquad a=2d^5,\qquad \phi=ax^2+\tfrac12(y^2+z^2).
\]
\end{lemma}
\begin{proof}
Since $A<\delta$ and $A+\delta\ge\frac32\delta>1$, the pair $(\frac12,1)$ applied to $F$ in $x_4$ has width $\max(A,1-A)<\delta$. By Lemmas~\ref{lem:normal} and~\ref{lem:lift}, failure forces $A=\frac12$ and $\lambda=0$. Then $t\ge A=\frac12$ and $t=d^5/a\le(m+1)d^4$. Hence
\[
\frac t\delta\le\frac{4(m+1)}{m+3}d^3<\frac{256(m+1)}{(m+3)^4}<1 .
\]
So the pair $(\frac12,1)$ applied to $G$ (target $\Delta=\delta$, as $K=1$) has width $t<\delta$, and failure forces $t=\frac12$ and $a=2d^5$. In particular $d^4\ge\frac1{2(m+1)}$. The coefficients of $x_4^2$ and $x_5^2$ in $\phi$ are $\rho_4=\frac12+ah^2$ and $\rho_5=\frac12+ak^2$. For $m=3$, $a<\frac d3<\frac29$ and $\rho_i<\frac59<8^{-1/4}\le d$. For $m=4$, $a<\frac d4<\frac17$ and $\rho_i<\frac{15}{28}<10^{-1/4}\le d$. So the pair $(\frac12,1)$ applied directly to $Q$ in $x_4$ and $x_5$ has width $\rho_i<d$. Failure forces $\rho_i=\frac12$, i.e.\ $h=k=0$.
\end{proof}

\begin{theorem}\label{thm:41}
The theorem holds when $(m,K)=(4,1)$.
\end{theorem}
\begin{proof}
Suppose not. By Lemma~\ref{lem:diag}, $a=2d^5$, $d^4\ge\frac1{10}$ and $d<\frac47$. Hence $a<2(\frac47)^5<\frac18$. The pair $(\frac12,2)$ applied to $Q$ in $x_3$ has width $\frac12-4a+\frac12=1-8d^5$. This is smaller than $d$ because $8d^5=4a\ge\frac45d$ and $\frac95d>1$ for $d\ge10^{-1/4}>\frac59$. The resonance $a=\frac18$ is excluded, so Lemma~\ref{lem:macbeath} gives $0<Q<d$, a contradiction.
\end{proof}

\begin{theorem}\label{thm:31}
The theorem holds when $(m,K)=(3,1)$.
\end{theorem}
\begin{proof}
Suppose not. By Lemma~\ref{lem:diag}, $\phi=ax^2+\frac12(y^2+z^2)$ with $a=2d^5$. From $a\ge d/4$ and $\delta=d+2a=d+4d^5<1$ we get $\frac{5946}{10000}<8^{-1/4}\le d<\frac{624}{1000}$. For each row of Table~\ref{tab:31int},
\[
\max\bigl(|H-2N^2\ell^5|,\,|H-2N^2u^5|\bigr)+\tfrac12<\ell .
\]
Since $d\mapsto2N^2d^5$ is increasing, for $\ell\le d\le u$ the distance $|H-2N^2d^5|$ is at most the larger of its two endpoint values. So the Macbeath width of $(H,N)$ in $x_3$ is less than $\ell\le d$ throughout $[\ell,u]$. So $0<Q<d$ unless $a=H/N^2$ is one of the five resonances listed.

\begin{table}[ht]
\centering
\begin{tabular}{cccc}
\toprule
$[\ell,u]$ & $N$ & $H$ & resonance $a=H/N^2$\\
\midrule
$[0.5946,\,0.5956]$ & $8$ & $19/2$ & none in range\\
$[0.5956,\,0.5964]$ & $10$ & $15$ & $3/20$\\
$[0.5964,\,0.6051]$ & $4$ & $5/2$ & $5/32$\\
$[0.6051,\,0.6166]$ & $3$ & $3/2$ & $1/6$\\
$[0.6166,\,0.6209]$ & $5$ & $9/2$ & $9/50$\\
$[0.6209,\,0.6240]$ & $4$ & $3$ & $3/16$\\
\bottomrule
\end{tabular}
\caption{Exact interval tests for $(3,1)$ (all endpoints are exact decimals).}\label{tab:31int}
\end{table}

At a resonance, fix $x_2=1$ and write $x_3=x+c_3$. The values are $j+\nu+\beta x-ax^2$ with $\beta\in H/N+\Lambda$, $\Lambda=\frac1N\Z+2a\Z$, a cyclic group with generator $g=\gcd(\frac1N,2a)$. For each of the finitely many $\beta$ in this coset, the residues of $\beta x-ax^2$ modulo $1$ are periodic in $x$. Table~\ref{tab:31gap} gives the exact largest circular gap over a full period. In the four rows with $(\text{gap})^5<a/2=d^5$, Lemma~\ref{lem:gapuse} gives $0<Q<d$.

\begin{table}[ht]
\centering
\begin{tabular}{cccc}
\toprule
$a$ & $g$ & largest gap & conclusion\\
\midrule
$3/20$ & $1/10$ & $1/4$ & solved\\
$5/32$ & $1/16$ & $1/4$ & solved\\
$9/50$ & $1/25$ & $4/25$ & solved\\
$3/16$ & $1/8$ & $1/2$ & solved\\
$1/6$ & $1/3$ & $2/3$ & see below\\
\bottomrule
\end{tabular}
\caption{Exact residue gaps at the $(3,1)$ resonances.}\label{tab:31gap}
\end{table}

It remains to treat $a=\frac16$, where $d^5=\frac1{12}$ and $d>\frac35$ (as $(\frac35)^5<\frac1{12}$). By the argument of Lemma~\ref{lem:normal} (row $(3/2,3)$), applied to $Q$ in $x_3$, failure for $x_2=\pm1$ and a shear of $x_3$ by $x_2$ give $(b,c):=(a_3,c_3)\in\{(\frac16,0),(0,\frac12)\}$. The Macbeath failures in $x_4$ and $x_5$ already forced $(e,c_y),(f,c_z)\in\{(0,\frac12),(\frac12,0)\}$, where $e=a_4$, $f=a_5$, $c_y=c_4$, $c_z=c_5$. For $x_2=n$, $x_3=x+c$, $x_4=y+c_y$, $x_5=z+c_z$, $\alpha=a_2$ and $\kappa=c_1$,
\begin{equation}\label{eq:31witness}
Q=nj+n^2\alpha+n\kappa+nb(x+c)-\tfrac16(x+c)^2+ne(y+c_y)-\tfrac12(y+c_y)^2+nf(z+c_z)-\tfrac12(z+c_z)^2 .
\end{equation}
Take $n=2$ first. Then~\eqref{eq:31witness} equals $2j+4\alpha+2\kappa+K(x,y,z)$, and $K$ is periodic modulo $2$ with periods $12$, $4$ and $4$ in $x$, $y$ and $z$. Enumerating one period shows that the residues of $K$ modulo $2$ have largest circular gap $\frac13$ in seven of the eight classes. Hence, by the argument of Lemma~\ref{lem:gapuse} on $\R/2\Z$, $0<Q<\frac35<d$ in those classes. The remaining class $b=e=f=0$, $c=c_y=c_z=\frac12$ has gap $\frac23$ at $n=2$. For it, \eqref{eq:31witness} is a fixed family with $|n|\le2$, $x,y,z\in\{-3,\dots,3\}$, $\mu=0$ and target $\frac35$, which has a dyadic cover. So $0<Q<\frac35<d$ in every class, a contradiction.
\end{proof}

\begin{remark}[Six variables]\label{rem:six}
The arguments of Section~\ref{sec:reduction} and of this section transfer to type $(1,5)$ at $d=(16|D|)^{1/6}$. There the negative block is quaternary, $\Phi=a(z_0+\mathbf h\cdot\mathbf z)^2+\Psi(\mathbf z)$ with $\Psi=b(z_1+\eta z_2+\zeta z_3)^2+c(z_2+\lambda z_3)^2+tz_3^2$. Lemma~\ref{lem:prune} with $(n,C)=(6,16)$ eliminates $m\ge4$ (as $7^5>2048\cdot5$), and with $k=2$ it leaves only $(m,K)\in\{(3,1),(2,1),(2,2),(1,1),(1,2)\}$. In the pair $(3,1)$ the analogue of Lemma~\ref{lem:diag} forces $b=c=t=\frac12$, $\eta=\zeta=\lambda=0$ and $a=2d^6$, so $\delta>\frac{31}{32}$ and $\frac a4<\frac1{18}$. The eight remaining shift classes are then settled without certificates. Seven are handled by residues modulo $2$ at $x_2=2$, as in Theorem~\ref{thm:31}. The class with all three shifts equal to $\frac12$ is handled by combining $x_2=\pm1$ with $x_2=2$. The details are in the supplementary note. The four remaining pairs are those where Sections~\ref{sec:third}--\ref{sec:22} need their computational tools. Since $\Gamma_{1,5}=16$ is known (see~\cite[\S1]{BKR}), we record this only as evidence that the reductions transfer between dimensions.
\end{remark}
\section{\texorpdfstring{The cases $(3,2)$, $(2,4)$ and $(2,3)$}{The cases (3,2), (2,4) and (2,3)}}\label{sec:third}

In this section and the next, $c_2\equiv0$, \eqref{eq:residual} holds and $\delta<1$. By rounding we may assume $\Delta<1$. Once $\Delta>t$ has been checked, $L$ is defined and rounding applied to $H$ lets us assume
\begin{equation}\label{eq:budget}
\Delta_L=\Delta+\tfrac{L^2-1}{4}t<1,\qquad\text{hence}\qquad \tfrac{L+3}{4}\Delta<1
\end{equation}
by~\eqref{eq:DeltaL}. We write $\mu=(a+A)/4$ for the constant of $G$. With $x_2=n$, $x_5=z+c_5$, $x_1=j+c_1$, $\alpha=a_2''$ and $\kappa=c_1$,
\begin{equation}\label{eq:Gfamily}
G=nj+n^2\alpha+n\kappa+ne(z+c)-q(z+c)^2-\mu,\qquad (q,e,c)=(t,a_5'',c_5).
\end{equation}
At a resonance of Lemma~\ref{lem:normal}, $q$ and $(e,c)$ take finitely many values, and~\eqref{eq:Gfamily} is a fixed family in the sense of Definition~\ref{def:cover}.

\begin{lemma}[Common scalar reduction]\label{lem:scalar}
Let $(m,K)\in\{(3,2),(2,4),(2,3),(2,2)\}$ and put
\[
C=C_{mK}=\frac{(m+3)(K+3)}{16},\qquad B=B_{mK}=\frac{2(m+1)(K+1)}{m+3}.
\]
Then $\Delta\ge Cd$ and $t\le Bd^3$. If $\Delta<1$, then $\Delta/t>1$. If moreover~\eqref{eq:budget} holds, then
\begin{equation}\label{eq:scalar}
\frac{C}{B(L+1)}\le d^2,\qquad d<\frac4{C(L+3)},
\end{equation}
and consequently $C^3(L+3)^2<16B(L+1)$. This leaves $L\le3$ for $(3,2)$ and $(2,4)$, and $L\le6$ for $(2,3)$. For $(2,2)$ also $\Delta\le5t$, so $L\le4$.
\begin{center}
\begin{tabular}{lccc}
\toprule
$(m,K)$ & $C_{mK}$ & $B_{mK}$ & possible $L$\\
\midrule
$(3,2)$ & $15/8$ & $4$ & $1,2,3$\\
$(2,4)$ & $35/16$ & $6$ & $1,2,3$\\
$(2,3)$ & $15/8$ & $24/5$ & $1,\dots,6$\\
$(2,2)$ & $25/16$ & $18/5$ & $1,2,3,4$\\
\bottomrule
\end{tabular}
\end{center}
\end{lemma}
\begin{proof}
We have $a\ge\frac d{m+1}$, so $\delta\ge d+\frac{(m-1)d}4=\frac{(m+3)d}4$. Also $A\ge\frac\delta{K+1}\ge\frac{(m+3)d}{4(K+1)}$ and $\Delta\ge\delta+\frac{(K-1)\delta}{4}=\frac{K+3}4\delta\ge Cd$. Hence $t=\frac{d^5}{2aA}\le Bd^3$. If $\Delta<1$ then $d<1/C$, so $\Delta/t\ge C/(Bd^2)>C^3/B>1$ (the last inequality is checked row by row). From $\Delta\le(L+1)t$ we get $Cd\le(L+1)Bd^3$, and from~\eqref{eq:budget} we get $Cd<\frac4{L+3}$. The two parts of~\eqref{eq:scalar} are compatible only if $C^3(L+3)^2<16B(L+1)$. This fails at $L=4$ for $(3,2)$ and $(2,4)$ and at $L=7$ for $(2,3)$. It then fails for all larger $L$, because $(L+3)^2/(L+1)$ increases. For $(2,2)$, $\delta\le3A$ gives $\Delta=\delta+\frac34A\le\frac{15}4A\le5t$.
\end{proof}

\begin{lemma}[Uniform gaps]\label{lem:U}
The exact values of $\gap_4(I)$ for the intervals below are as stated.
\begin{center}
\begin{tabular}{lccc}
\toprule
Row & $I$ & $\gap_4(I)$ & used in\\
\midrule
U1 & $[7/50,\,157/1000]$ & $57/100$ & $(2,3)$, $L=1$, $F$\\
U2 & $[22/125,\,1/5]$ & $3/5$ & $(2,3)$, $L=1$, $F$\\
U3 & $[73/500,\,63/400]$ & $23/40$ & $(2,3)$, $L=3$\\
U4 & $[22/125,\,27/125]$ & $3/5$ & $(2,3)$, $L=3$\\
U5 & $\{1/5\}$ & $3/5$ & $(2,2)$, $L=3$\\
\bottomrule
\end{tabular}
\end{center}
\end{lemma}
\begin{proof}
Exact evaluation of Proposition~\ref{prop:gapexact} in rational arithmetic.
\end{proof}

\begin{lemma}[Shared resonant covers]\label{lem:R}
Each family~\eqref{eq:Gfamily} in Table~\ref{tab:R}, with $|n|\le6$, $|z|\le10$, $\mu\in[\mu_-,\mu_+]$ and the single target $(0,W)$, has a dyadic cover. Hence for all $\alpha,\kappa\in\R$ and $\mu\in[\mu_-,\mu_+]$ there are $n\ne0$, $z$, $j$ with $0<G<W$.
\end{lemma}

\begin{table}[ht]
\centering
\begin{tabular}{lccccc}\toprule Row & $q$ & $(e,c)$ & $\mu_-$ & $\mu_+$ & $W$\\\midrule
R0a & $1/2$ & $(0,1/2)$ & $3/32$ & $3/32$ & $15/16$\\
R0b & $1/2$ & $(1/2,0)$ & $3/32$ & $3/32$ & $15/16$\\
R1 & $1/2$ & $(0,1/2)$ & $301/4800$ & $27/280$ & $19/20$\\
R2 & $1/2$ & $(1/2,0)$ & $301/4800$ & $27/280$ & $19/20$\\
R3 & $1/2$ & $(0,1/2)$ & $1457/19200$ & $23/180$ & $22/25$\\
R4 & $1/6$ & $(1/6,0)$ & $1/25$ & $53/1000$ & $3/5$\\
R5 & $1/6$ & $(0,1/2)$ & $1/25$ & $53/1000$ & $3/5$\\
R6 & $1/4$ & $(0,0)$ & $217/3600$ & $23/225$ & $7/10$\\
R7 & $1/4$ & $(1/4,1/2)$ & $217/3600$ & $23/225$ & $7/10$\\
R8 & $1/6$ & $(1/6,0)$ & $5053/96000$ & $351/5000$ & $61/100$\\
R9 & $1/6$ & $(0,1/2)$ & $5053/96000$ & $351/5000$ & $61/100$\\
\bottomrule\end{tabular}
\caption{Shared resonant families, used in $(3,2)$, $(2,4)$ and $(2,3)$.}\label{tab:R}
\end{table}

\begin{theorem}\label{thm:32}
The theorem holds when $(m,K)=(3,2)$.
\end{theorem}
\begin{proof}
Here $\frac d4\le a<\frac d3$, $\delta=d+2a$, $\frac\delta3\le A<\frac\delta2$ and $\Delta=\delta+\frac34A$. Lemma~\ref{lem:scalar} gives $\Delta\ge\frac{15}8d$ and $t\le4d^3$.

If $\delta+A>1$, the pair $(\frac12,1)$ applied to $F$ in $x_4$ has width $1-A<\delta$ (note $A<\frac\delta2<\frac12$), and the resonance $A=\frac12$ is excluded. So assume $\delta+A\le1$. Since $\delta\ge\frac32d$ and $A\ge\frac\delta3\ge\frac d2$, this gives $d\le\frac12$. By Lemma~\ref{lem:scalar}, $L\in\{1,2,3\}$, $\frac{15}{32(L+1)}\le d^2$ and $d<\frac{32}{15(L+3)}$.

\emph{$L=1$}, $\sqrt{15/64}\le d\le\frac12$. Then $\Delta\ge\frac{15}{8}\sqrt{15/64}>\frac23$, $\frac\Delta2\le t<\Delta$ and $t\le4d^3\le\frac12$. The pair $(\frac12,1)$ applied to $G$ has width $\max(t,1-t)<\Delta$. At the resonance $t=\frac12$ every inequality in $t\le4d^3\le\frac12$ is an equality, so $(d,a,A)=(\frac12,\frac18,\frac14)$, $\Delta=\frac{15}{16}$ and $\mu=\frac3{32}$. By Lemma~\ref{lem:normal}, $(e,c)\in\{(0,\frac12),(\frac12,0)\}$, and rows R0a, R0b give $0<G<\frac{15}{16}=\Delta$.

\emph{$L=2$}, $\sqrt{5/32}\le d<\frac{32}{75}$. Here $\Delta\ge\frac{15}8\sqrt{5/32}>\frac{37}{50}$. Take the pair $(1,2)$ for $G$. If $4t\le1$, its width is $\frac32-4t\le\frac32-\frac43\Delta<\Delta$, because $\Delta>\frac9{14}$. If $4t\ge1$, its width is $4t-\frac12\le16d^3-\frac12<\frac{15}{8}d\le\Delta$. Indeed $16d^3-\frac{15}8d$ increases for $d\ge\sqrt{5/32}$, and $(16(\frac{32}{75})^2-\frac{15}8)\frac{32}{75}<\frac12$. Failure forces $t=\frac14$ and $(e,c)\in\{(0,0),(\frac14,\frac12)\}$. Then $aA=2d^5$, so $A\le8d^4$ and
\[
\tfrac{217}{3600}<\tfrac{117}{1600}\le\tfrac{3d}{16}\le\mu\le\tfrac14\bigl(\tfrac d3+8d^4\bigr)\le\tfrac{3222152}{31640625}<\tfrac{23}{225}.
\]
Rows R6, R7 give $0<G<\frac7{10}<\Delta$.

\emph{$L=3$}, $\sqrt{15/128}\le d<\frac{16}{45}$. Here $\Delta\ge\frac{15}{8}\sqrt{15/128}>\frac{16}{25}$. Take the pair $(\frac32,3)$. If $9t\le\frac32$, the width is $2-9t\le2-\frac94\Delta<\Delta$, because $\Delta>\frac8{13}$. If $9t\ge\frac32$, it is $9t-1\le36d^3-1<36(\frac{16}{45})^3-1<\frac{15}{8}\cdot\frac{34}{100}<\Delta$. Failure forces $t=\frac16$ and $(e,c)\in\{(\frac16,0),(0,\frac12)\}$, and then $aA=3d^5$. The function $a\mapsto a+3d^5/a$ is convex on $[\frac d4,\frac d3]$, so
\[
\tfrac{5053}{96000}<\tfrac{51}{800}\le\tfrac{3d}{16}\le\mu\le\tfrac14\max\bigl(\tfrac d4+12d^4,\ \tfrac d3+9d^4\bigr)\le\tfrac{95911}{1366875}<\tfrac{351}{5000},
\]
where the right side is evaluated at $d=\frac{16}{45}$ (both expressions increase with $d$). Rows R8, R9 give $0<G<\frac{61}{100}<\Delta$.
\end{proof}

\begin{theorem}\label{thm:24}
The theorem holds when $(m,K)=(2,4)$.
\end{theorem}
\begin{proof}
Here $\frac d3\le a<\frac d2$, $\delta=d+\frac34a$, $\frac\delta5\le A<\frac\delta4$ and $\Delta=\delta+\frac{15}4A$. By Lemma~\ref{lem:scalar}, $\Delta\ge\frac{35}{16}d$, $t\le6d^3$, $L\in\{1,2,3\}$, $\frac{35}{96(L+1)}\le d^2$ and $d<\frac{64}{35(L+3)}$.

\emph{$L=1$.} Here $\Delta\ge\frac{35}{16}\sqrt{35/192}>\frac23$ and $\frac\Delta2\le t<\Delta$, so $(\frac12,1)$ has width $\max(t,1-t)<\Delta$. At the resonance $t=\frac12$, $(e,c)\in\{(0,\frac12),(\frac12,0)\}$ and $aA=d^5\ge\frac{d^2}{12}$. Hence $d>\frac{43}{100}$ and $\Delta\ge\frac{35}{16}\,12^{-1/3}>\frac{19}{20}$, as $(\frac{19}{20}\cdot\frac{16}{35})^3<\frac1{12}$. Also $\mu\in[\frac{7d}{48},\frac{27d}{128}]\subset[\frac{301}{4800},\frac{27}{280}]$, using $A<\frac\delta4<\frac{11d}{32}$. Rows R1, R2 give $0<G<\frac{19}{20}<\Delta$.

\emph{$L=2$.} Here $d\ge\sqrt{35/288}>\frac{34}{100}$, so $t\ge\frac\Delta3\ge\frac{35d}{48}\ge\frac{35}{48}\sqrt{35/288}>\frac14$. The pair $(1,2)$ has width $4t-\frac12\le24d^3-\frac12\le24(\frac{64}{175})^3-\frac12<\frac{35}{16}\cdot\frac{34}{100}<\Delta$, and the resonance $t=\frac14$ is excluded. No certificate is needed.

\emph{$L=3$.} Here $\sqrt{35/384}\le d<\frac{32}{105}$ and $\Delta>\frac8{13}$. For the pair $(\frac32,3)$, the width is $2-9t\le2-\frac94\Delta<\Delta$ if $9t\le\frac32$. Otherwise it is $9t-1\le54(\frac{32}{105})^3-1<\frac{35}{16}\cdot\frac3{10}<\Delta$. At the resonance $t=\frac16$, $(e,c)\in\{(\frac16,0),(0,\frac12)\}$ and $aA=3d^5\ge\frac{d^2}{12}$, so $d>\frac3{10}$ and $\Delta\ge\frac{35}{16}36^{-1/3}>\frac{33}{50}$. Moreover
\[
\mu\in\bigl[\tfrac{7d}{48},\tfrac{27d}{128}\bigr]\subset\bigl[\tfrac7{160},\tfrac9{140}\bigr]\subset\bigl[\tfrac1{25},\tfrac{53}{1000}\bigr]\cup\bigl[\tfrac{5053}{96000},\tfrac{351}{5000}\bigr].
\]
Rows R4, R5 (target $\frac35$) and R8, R9 (target $\frac{61}{100}$) give $0<G<\Delta$.
\end{proof}

\begin{theorem}\label{thm:23}
The theorem holds when $(m,K)=(2,3)$.
\end{theorem}
\begin{proof}
Here $\frac d3\le a<\frac d2$, $\delta=d+\frac34a\ge\frac54d$, $\frac\delta4\le A<\frac\delta3$ and $\Delta=\delta+2A$. By Lemma~\ref{lem:scalar}, $\Delta\ge\frac{15}8d$, $t\le\frac{24}5d^3$, $1\le L\le6$ and
\begin{equation}\label{eq:23d}
\frac{25}{64(L+1)}\le d^2,\qquad d<\frac{32}{15(L+3)}.
\end{equation}

\emph{$L=1$.} Here $\Delta\ge\frac{15}8\sqrt{25/128}>\frac23$ and $\frac\Delta2\le t<\Delta$, so $(\frac12,1)$ has width less than $\Delta$, and failure forces $t=\frac12$ with $(e,c)\in\{(0,\frac12),(\frac12,0)\}$. Then $aA=d^5\ge\frac{5d^2}{48}$, so $d>\frac{47}{100}$, $\Delta\ge\frac{15}8d>\frac{22}{25}$ and $\mu\in[\frac{31d}{192},\frac{23d}{96}]\subset[\frac{1457}{19200},\frac{23}{180}]$. For $(e,c)=(0,\frac12)$, row R3 gives $0<G<\frac{22}{25}<\Delta$.

For $(e,c)=(\frac12,0)$ we work with $F$ at $x_2=1$, $x_5=0$. Its values are $j+\nu+\beta y-Ay^2$, and $\frac7{50}<\frac{5d}{16}\le A<\frac\Delta5<\frac15$. On $[\frac7{50},\frac{157}{1000}]$, row U1 gives $0<F<\frac{57}{100}<\delta$. On $[\frac{22}{125},\frac15]$, we have $d^4=aA/d\ge\frac A3\ge\frac{22}{375}>(\frac{12}{25})^4$, so $\delta\ge\frac54d>\frac35$ and row U2 applies. On $[\frac{157}{1000},\frac{22}{125}]$, the pair $(\frac32,3)$ has width at most $\frac{587}{1000}<\frac54\cdot\frac{47}{100}<\delta$. By Lemma~\ref{lem:normal} for $F$ (with $x_5\in\{0,1\}$), failure forces $A=\frac16$, $\lambda=0$ and $(a_4',c_4)\in\{(\frac16,0),(0,\frac12)\}$. Then $a=6d^5$, $\frac{97}{200}<d<\frac{511}{1000}$, $\delta=d+\frac92d^5>\frac35$, and $F$ is the family~\eqref{eq:Gfamily} with $q=\frac16$ and $\mu=\frac a4=\frac32d^5\in[\frac1{25},\frac{53}{1000}]$. Rows R4, R5 give $0<F<\frac35<\delta$.

\emph{$L=2$.} By~\eqref{eq:23d}, $d^2\ge\frac{25}{192}$ and $d<\frac{32}{75}$, so $\Delta>\frac9{14}$ and $t<\frac\Delta2<\frac12$. If $t\le\frac14$, the pair $(1,2)$ has width $\frac32-4t\le\frac32-\frac43\Delta<\Delta$. If $t\ge\frac14$, then $d^3\ge\frac5{96}$, so $d>\frac{28}{75}$ and $\Delta>\frac7{10}$. Then either $\Delta+t>1$, and $(\frac12,1)$ has width $1-t<\Delta$ with $t=\frac12$ excluded, or $t\le1-\Delta$, and $(1,2)$ has width $4t-\frac12\le\frac72-4\Delta<\Delta$. Failure forces $t=\frac14$ and $(e,c)\in\{(0,0),(\frac14,\frac12)\}$, with $\mu\in[\frac{31d}{192},\frac{23d}{96}]\subset[\frac{217}{3600},\frac{23}{225}]$. Rows R6, R7 give $0<G<\frac7{10}<\Delta$.

\emph{$L=3$.} Here $\frac5{16}\le d<\frac{16}{45}$, $\Delta\ge\frac{75}{128}$ and $\frac{73}{500}<\frac\Delta4\le t\le\frac{24}5d^3<\frac{27}{125}$. For $t\le\frac{63}{400}$, row U3 gives gap $\frac{23}{40}<\frac{75}{128}\le\Delta$. For $t\ge\frac{22}{125}$, $d^3\ge\frac{11}{300}>(\frac8{25})^3$, so $\Delta>\frac35$ and row U4 applies. On $[\frac{63}{400},\frac{22}{125}]$, the pair $(\frac32,3)$ has width at most $\max(\frac{233}{400},\frac{73}{125})<\frac{75}{128}$. At the resonance $t=\frac16$, $(e,c)\in\{(\frac16,0),(0,\frac12)\}$ and $aA=3d^5$. Then $d^3\ge\frac5{144}$, $d>\frac{163}{500}$ and $\Delta>\frac{61}{100}$. By convexity of $a\mapsto a+3d^5/a$ on $[\frac d3,\frac d2]$,
\[
\tfrac{5053}{96000}\le\tfrac{31d}{192}\le\mu\le\tfrac14\max\bigl(\tfrac d3+9d^4,\ \tfrac d2+6d^4\bigr)<\tfrac7{100}\qquad(d<\tfrac{16}{45}).
\]
Rows R8, R9 give $0<G<\frac{61}{100}<\Delta$.

\emph{$L=4,5,6$.} By $t\ge\Delta/(L+1)$ and~\eqref{eq:23d},
\[
\Delta+t\ge\frac{L+2}{L+1}\cdot\frac{15d}8\ge\frac{75(L+2)}{64(L+1)^{3/2}}\ge\frac{75}{56\sqrt7}>\frac12 .
\]
The middle bound decreases in $L$ because $(x+1)/x^{3/2}$ decreases for $x>0$; the last inequality is $75>28\sqrt7$, i.e.\ $75^2=5625>5488=28^2\cdot7$. Lemma~\ref{lem:highq} gives $0<G<\Delta$.

In every branch Lemma~\ref{lem:lift} lifts the value of $G$ (or $F$) to $0<Q<d$.
\end{proof}
\section{\texorpdfstring{The case $(2,2)$}{The case (2,2)}}\label{sec:22}

\subsection{The domain}
Let $c_2\equiv0$, \eqref{eq:residual}, $2<d/a\le3$ and $2<\delta/A\le3$, where now $\delta=d+\frac34a$ and
\[
\Delta=d+\tfrac34(a+A)=d+3\mu,\qquad \mu:=\tfrac{a+A}4 .
\]
We may assume $\delta<1$, $\Delta<1$ and, once $L$ is defined,~\eqref{eq:budget}. The reduction inequalities give
\begin{equation}\label{eq:22dom}
\tfrac d3\le a<\tfrac d2,\qquad \tfrac\delta3\le A<\tfrac\delta2,\qquad A\ge\tfrac34a,\qquad t\ge\tfrac34A,\qquad A\ge\tfrac{5d}{12},
\end{equation}
and Lemma~\ref{lem:scalar} gives $\Delta\ge\frac{25d}{16}$, $t\le\frac{18}5d^3$, $L\in\{1,2,3,4\}$ and
\begin{equation}\label{eq:22d}
\sqrt{\tfrac{125}{288(L+1)}}\le d<\tfrac{64}{25(L+3)} .
\end{equation}
We shall use the following lower bounds for the target $\Delta$ of $G$; each is valid on the whole $(2,2)$ domain:
\begin{equation}\label{eq:22targets}
\Delta>Lt,\qquad \Delta\ge d_0+3\mu\ \text{ for any } d_0\le d,\qquad \Delta>\tfrac{121}{19}\mu .
\end{equation}
For the last one, $4\mu=a+A<\frac d2+\frac\delta2<\frac d2+\frac12\cdot\frac{11d}{8}=\frac{19}{16}d$, so $\Delta=d+3\mu>\frac{64}{19}\mu+3\mu$. Finally $t\le\frac{18}5d^3$ gives $d\ge(\frac{5t}{18})^{1/3}$; we write $d_0(t)$ for a rational lower bound for $(\frac{5t}{18})^{1/3}$.

\subsection{\texorpdfstring{The branch $L=1$}{The branch L=1}}
Here $\frac\Delta2\le t<\Delta$, $\Delta+t\ge\frac{75}{32}d>1$ and $d\ge\sqrt{125/576}>\frac{32}{75}$, so the pair $(\frac12,1)$ applied to $G$ has width $\max(t,1-t)<\Delta$. Failure forces $t=\frac12$ and $(a_5'',c_5)\in\{(0,\frac12),(\frac12,0)\}$. Then $d^5=aA$, so $A\le3d^4$, $d\ge(\frac5{36})^{1/3}>\frac{18}{35}$ and $\delta\ge\frac54d$.

If $A\le\frac14$, then $\delta+4A\ge\frac{35}{12}d>\frac32$. If $\frac14\le A\le\frac3{10}$, then $\delta\ge\frac54(A/3)^{1/4}>4A-\frac12$, because $g(A)=\frac54(A/3)^{1/4}-4A+\frac12$ is concave and positive at both endpoints (at $\frac3{10}$ use $(\frac7{10})^4<(\frac54)^4/10$). In both cases the pair $(1,2)$ applied to $F$ has width $|1-4A|+\frac12<\delta$. If $A\ge\frac3{10}$, then $A+\delta\ge A+\frac54(A/3)^{1/4}>1$ and $\delta>2A$, so the pair $(\frac12,1)$ has width $\max(A,1-A)<\delta$. Its resonance $A=\frac12$ is impossible, since $\delta+\frac34A<1$ and $\delta>2A$ give $A<\frac4{11}$.

Thus failure forces $A=\frac14$. By Lemma~\ref{lem:normal}, $\lambda=0$ and $(a_4',c_4)\in\{(0,0),(\frac14,\frac12)\}$, and $(a_5',c_5)=(a_5'',c_5)$ is as above. Now $a=4d^5$, $\frac16<a<\frac3{11}$, $d\ge12^{-1/4}$ and $\delta\ge\frac54\,12^{-1/4}>\frac23$ (as $(\frac8{15})^4<\frac1{12}$). With $x_2=n$, $x_4=y+c_4$ and $x_5=z+c_5$,
\[
F=nj+n^2\alpha+n\kappa+ne_4(y+c_4)+ne_5(z+c_5)-\tfrac14(y+c_4)^2-\tfrac12(z+c_5)^2-\tfrac a4 .
\]
For each of the four choices of $(e_4,c_4,e_5,c_5)$ this is a fixed family with $|n|\le4$, $|y|,|z|\le3$, $\mu=\frac a4\in[\frac1{24},\frac3{44}]$ and target $\frac23$. The four dyadic covers give $0<F<\frac23<\delta$.

\subsection{\texorpdfstring{The branch $L=4$}{The branch L=4}}
Here $t\ge\frac\Delta5$ and, by~\eqref{eq:22d}, $d^2\ge\frac{25}{288}$. Hence
\[
\Delta+t\ge\tfrac65\Delta\ge\tfrac{15}8d\ge\frac{25}{32\sqrt2}>\frac12 .
\]
Lemma~\ref{lem:highq} gives $0<G<\Delta$.

\subsection{\texorpdfstring{The branches $L=2,3$}{The branches L=2,3}}
\emph{Partition.} By~\eqref{eq:22d}, $t\ge\Delta/(L+1)\ge\frac{25d}{16(L+1)}$ and $t<\Delta/L<\frac4{L(L+3)}$. So the parameter $t$ lies in $[\frac{3803}{19200},\frac25]$ for $L=2$ and in $[\frac{1647}{12800},\frac29]$ for $L=3$; these endpoints are $\frac{25\ell_L}{16(L+1)}$ with $\ell_2=\frac{3803}{10^4}$ and $\ell_3=\frac{3294}{10^4}$. These ranges are partitioned into $8$ and $15$ rational intervals $[t_-,t_+]$. On each interval put
\[
\delta_-=\max\bigl(Lt_-,\ \tfrac{25}{16}\ell_L,\ \tfrac{25}{16}d_0(t_-)\bigr)\le\Delta .
\]
Then either $\gap_M([t_-,t_+])<\delta_-$ with $M\in\{4,6\}$ ($5$ and $11$ intervals), or some Macbeath pair $(H,N)$ has width less than $\delta_-$ at both endpoints (the width is convex in $t$). Only four intervals contain the resonance of their pair:
\begin{center}
\begin{tabular}{cccc}
\toprule
$L$ & interval & pair & resonance\\
\midrule
$2$ & $[3803/19200,\ 82423/409600]$ & $(5,5)$ & $t=1/5$\\
$2$ & $[6363/25600,\ 8411/30720]$ & $(1,2)$ & $t=1/4$\\
$3$ & $[10061/61440,\ 312607/1843200]$ & $(3/2,3)$ & $t=1/6$\\
$3$ & $[30541/153600,\ 743761/3686400]$ & $(5,5)$ & $t=1/5$\\
\bottomrule
\end{tabular}
\end{center}
The complete list of intervals, with exact endpoints and bounds, is part of the verification data. By Lemma~\ref{lem:normal}, $0<G<\Delta$ except at these four values of $t$ with the normal forms of Table~\ref{tab:normal}.

\emph{The resonance $t=\frac15$ for $L=3$.} Here $\Delta>3t=\frac35$, and row U5 gives $\gap_4(\{\frac15\})=\frac35<\Delta$. No certificate is needed.

\emph{The other resonances.} At $t=t_0\in\{\frac15,\frac14\}$ ($L=2$) and $t_0=\frac16$ ($L=3$), $G$ is the family~\eqref{eq:Gfamily} with $q=t_0$ and $(e,c)$ a normal form. We use $d\ge d_0:=\max(\ell_L,d_0(t_0))$, $\mu\ge\frac{3d_0}{16}$ and $\mu\le\frac13(\Delta_{\max}-d_0)$, where $\Delta_{\max}=\min((L+1)t_0,\,1-\frac{L^2-1}4t_0)$, together with the three targets~\eqref{eq:22targets}: $(c_1,c_0)=(3,d_0)$, $(\frac{121}{19},0)$ and $(0,Lt_0)$. Table~\ref{tab:22res} lists the seven families.

\begin{table}[ht]
\centering\small
\begin{tabular}{lcccccc}\toprule $L$ & $t$ & $(e,c)$ & $d_0$ & $\mu_-$ & $\mu_+$ & pinch leaves\\\midrule
2 & $1/4$ & $(0,0)$ & $41103/100000$ & $123309/1600000$ & $11299/100000$ & no cover needed\\
2 & $1/4$ & $(1/4,1/2)$ & $41103/100000$ & $123309/1600000$ & $11299/100000$ & --\\
2 & $1/5$ & $(0,0)$ & $38157/100000$ & $114471/1600000$ & $7281/100000$ & --\\
2 & $1/5$ & $(1/10,1/4)$ & $38157/100000$ & $114471/1600000$ & $7281/100000$ & --\\
2 & $1/5$ & $(1/5,1/2)$ & $38157/100000$ & $114471/1600000$ & $7281/100000$ & --\\
3 & $1/6$ & $(0,1/2)$ & $35907/100000$ & $107721/1600000$ & $92279/900000$ & --\\
3 & $1/6$ & $(1/6,0)$ & $35907/100000$ & $107721/1600000$ & $92279/900000$ & pinch ball\\
\bottomrule\end{tabular}
\caption{Resonant families for $(2,2)$, $L=2,3$. For $(2,\frac14,0,0)$ no cover is needed (Theorem~\ref{thm:fail}); for $(3,\frac16,\frac16,0)$ the cover has pinch leaves in the ball $\bar B((\frac7{12},\frac12),\frac12)$.}\label{tab:22res}
\end{table}

\begin{theorem}[Failure sets]\label{thm:fail}
In the resonant families of Table~\ref{tab:22res}, write $(\alpha,\kappa)=(a_2'',c_1)$ modulo $\Z^2$. Then $G$ has a value in $(0,\Delta)$ at every admissible point, except exactly at the following points:
\begin{enumerate}
\item for $(L,t_0,e,c)=(2,\frac14,0,0)$: $(\alpha,\kappa)\in\{(0,0),(\frac12,\frac12)\}$ and $d+a+A\le\frac34$;
\item for $(3,\frac16,\frac16,0)$: $(\alpha,\kappa)=(\frac7{12},\frac12)$, $\mu\ge\frac1{12}$ and $d+a+A\le\frac34$.
\end{enumerate}
In the other five families, $G$ has a value in $(0,\Delta)$ at every admissible point.
\end{theorem}
\begin{proof}
\emph{Values at the three points.} At $(\alpha,\kappa)=(0,0)$ or $(\frac12,\frac12)$ in the first family, $G+\mu=nj+n^2\alpha+n\kappa-\frac14z^2\in\Z-\{0,\frac14\}$, since $n(n+1)/2\in\Z$ and $z^2\equiv0,1\pmod4$. So the values of $G+\mu$ in $(0,1]$ are $\frac34$ and $1$, and the window $(\mu,\mu+\Delta)=(\mu,d+a+A)$ misses them iff $d+a+A\le\frac34$. At $(\frac7{12},\frac12)$ in the second family, put $R=12(G+\mu)=12nj+7n^2+6n+2nz-2z^2$. For odd $n$, $R\equiv1$ or $9\pmod{12}$. For even $n=2k$, $R\equiv0,4,6$ or $10\pmod{12}$, and $R\equiv4kz-2z^2\pmod8$, which is never $4\pmod 8$. If $R\equiv6\pmod8$ then $k$ is even and $z$ odd, and then $R\equiv14\pmod{16}$. Hence the values of $G+\mu$ in $(0,1]$ are exactly $\frac1{12},\frac34,\frac56$, attained at $(n,z,j)=(1,0,-1),(1,2,0),(2,3,-1)$. The window misses them iff $\mu\ge\frac1{12}$ and $d+4\mu\le\frac34$.

\emph{Ladders.} We apply Lemma~\ref{lem:ladder} with the lower bound $\Delta\ge d_0+3\mu$ on the $\mu$-range of Table~\ref{tab:22res}.
\begin{itemize}
\item At $(0,0)$: ascending rungs $(0,n)$ for all $n\ge1$ (take $z=j=0$), so $\rho_k\le2$. Descending rungs $(\frac34,\frac{z^2+3}4)$ for odd $z$, with $j=1$ (and $j=-1$ for $-n$); here $\rho_k\le3$.
\item At $(\frac12,\frac12)$: ascending rungs $(0,n)$ for $n=1,2$ and all odd $n$, so $\rho_k\le2$. For odd $n$ the term $n(n+1)/2$ lies in $n\Z$. The descending rungs are those of $(0,0)$, since $n=\frac{z^2+3}4$ is odd.
\item At $(\frac7{12},\frac12)$: $V_{-n}=V_n$, because $R(-n,-z,-j-1)=R(n,z,j)$. Ascending rungs $(\frac1{12},n)$ for $n=1,3$ and $n=2z^2+1$ with $z\equiv0,2\pmod3$ and $j=-1-\frac{7z^2+z}6$; the ratios are at most $3$. Descending rungs $(\frac34,1),(\frac56,2),(\frac34,3),(\frac56,6),(\frac34,9),(\frac34,11),(\frac34,17),(\frac34,19)$, followed by $(\frac34,6w^2+3)$ for $w\ge2$ with $z=3w$ and $j=-\frac{7w^2+w+4}2$. The ratios in the tail are at most $\frac{57}{27}$ and decrease.
\end{itemize}
The handoff inequalities (A) and (D) hold for every rung pair at both endpoints of the $\mu$-range. The radii are $r\ge\frac{63701}{10^5}$ at the first two points and $r\ge\frac{654047}{1200000}$ at the third, all larger than $\frac12$. So Lemma~\ref{lem:ladder} applies on the whole \emph{closed} $\ell^1$-ball of radius $\frac12$ about each point, with the centre removed. For the family $(2,\frac14,0,0)$ the two closed balls of radius $\frac12$ about $(0,0)$ and $(\frac12,\frac12)$ cover the torus: in the coordinates $(\alpha+\kappa,\alpha-\kappa)$ the centres form $\Z^2$, the $\ell^1$-norm becomes the $\ell^\infty$-norm, and closed $\ell^\infty$-balls of radius $\frac12$ about $\Z^2$ cover $\R^2$. (Open balls would not suffice: $(\frac12,0)$ is at distance exactly $\frac12$ from both centres.) For the other six families, the dyadic covers (with pinch leaves for $(3,\frac16,\frac16,0)$) and Lemma~\ref{lem:cover} give the rest.
\end{proof}

\emph{The pinch points.} At the three points of Theorem~\ref{thm:fail} we return to $F$ and keep the coupling between its parameters. Write
\[
F=(x_1+a_2''x_2+e\,x_5)x_2-A(x_4-qx_2+\lambda x_5)^2-t_0x_5^2-\tfrac a4,\qquad q=\frac{a_4'}{2A},\ e=a_5'' .
\]
This identity follows by completing the square in $x_4$. Put $x_2=n$, $x_4=y+c_4$, $x_5=z+c$ and $C=c_4+\lambda c$. The substitutions $y\mapsto y+rn$ and $y\mapsto y+r'$ ($r,r'\in\Z$) leave $a_2''$, $e$, $c_1$ and the $G$-parameters unchanged, so we may assume $q,C\in[-\frac12,\frac12]$. With the moment coordinates
\[
p=A\lambda,\ v=A\lambda^2,\ r=Aq,\ s=A\lambda q,\ u=AC,\ w=A\lambda C,\ m=AC^2,\ Q_2=Aq^2,\ X=AqC,
\]
$\alpha_G=a_2''$ and $\kappa_G=c_1$, we obtain
\begin{multline}\label{eq:corner}
F=nj+n^2(\alpha_G-Q_2)+n(\kappa_G+ec+2X)+enz-Ay^2-2pyz+2rny-2uy\\
-vz^2+2snz-2wz-m-t_0(z+c)^2-\tfrac a4 .
\end{multline}
This is affine in $\theta=(d,a,A,P,p,v,r,s,u,w,m,Q_2,X,\alpha_G,\kappa_G)$, where $P=aA$. Appendix~\ref{app:corner} lists rational inequalities valid at every failure point of Theorem~\ref{thm:fail}. For each of the three pinch points there is a Farkas tree for these inequalities and the target $\delta=d+\frac34a$ ($713$ nodes in total). By Lemma~\ref{lem:farkas} we get $0<F<\delta$.

\begin{theorem}\label{thm:22}
The theorem holds when $(m,K)=(2,2)$.
\end{theorem}
\begin{proof}
Combine the branches above with Lemma~\ref{lem:lift}.
\end{proof}
\section{\texorpdfstring{The case $(2,1)$}{The case (2,1)}}\label{sec:21}

Let $c_2\equiv0$, \eqref{eq:residual}, $\frac d3\le a<\frac d2$ and $\frac\delta2\le A<\delta$ with $\delta=d+\frac34a$. Then $d\le\frac34$, $\delta<1$ and
\begin{equation}\label{eq:21bounds}
A\ge\tfrac{5d}8,\qquad A^2\le\tfrac{2d^5}{3a}\le2d^4,\qquad \tfrac{5\sqrt2}{16}\le d,\qquad \tfrac{15d}{32}\le t\le\tfrac{12}5d^3 .
\end{equation}

\subsection{\texorpdfstring{Low $A$}{Low A}}
If $A+\delta\le1$, then $\frac{15}8d\le A+\delta\le1$, so $\frac{5\sqrt2}{16}\le d\le\frac8{15}$. With $x_2=1$, expanding the shifts of the negative coordinates gives
\begin{equation}\label{eq:shiftexp}
Q=j+\nu+\beta\cdot w-\phi(w),\qquad w\in\Z^3,
\end{equation}
with $\beta\in\R^3$ and $\nu\in\R$ arbitrary; we may take $\beta\in[-\frac12,\frac12]^3$ and $\nu\in[0,1]$.

\begin{lemma}\label{lem:21ternary}
If moreover $A+\delta\le1$, then for all $\beta,\nu$ there are $j\in\Z$ and $w\in\{-2,\dots,2\}^3$ with $0<j+\nu+\beta\cdot w-\phi(w)<d$.
\end{lemma}
\begin{proof}
Two Farkas trees of component~I prove this on the ternary outer polytopes of Appendix~\ref{app:ternary}. The two trees cover the sign regions $k\ge0$ and $k\le0$, after the reflection that makes $h\ge0$.
\end{proof}

\subsection{\texorpdfstring{High $A$}{High A}}
Assume $A+\delta>1$, and suppose that no $X\equiv c$ has $0<Q(X)<d$. The pair $(\frac12,1)$ applied to $F$ in $x_4$ has width $\max(A,1-A)<\delta$. By Lemma~\ref{lem:normal}, failure forces $A=\frac12$ and $\lambda=0$. Then $t\ge\frac12$ and $t=d^5/a\le3d^4$, so
\[
d\ge6^{-1/4},\qquad \tfrac12\le t\le\tfrac{243}{256}<1,\qquad \delta\ge\tfrac54\,6^{-1/4}>\tfrac34 .
\]
For $G$ (with $K=1$, target $\delta$), take the nearest $H\in\{\frac12,1\}$ to $t$. The width is at most $\frac34<\delta$, so failure forces $t=H$, and since $t<1$, $t=\frac12$. Hence $a=2d^5$. The coefficient of $x_4^2$ in $\phi$ is $\rho_4=\frac12+ah^2\le\frac12+\frac a4<\frac12+\frac d8<d$, because $d\ge6^{-1/4}>\frac47$. So the pair $(\frac12,1)$ applied to $Q$ in $x_4$ forces $h=0$, and likewise $k=0$. Using both signs of $x_2$,
\begin{equation}\label{eq:21diag}
\phi=ax_3^2+\tfrac12(x_4^2+x_5^2),\quad a=2d^5,\qquad (a_i,c_i)\in\{(0,\tfrac12),(\tfrac12,0)\}\ (i=4,5).
\end{equation}

\begin{lemma}\label{lem:21diag}
Let $a=2d^5$, $\frac d3\le a\le\frac d2$, $a+d\le1$ and $d\ge6^{-1/4}$. For all real $a_2,a_3,c_1,c_3$ and all pairs~\eqref{eq:21diag}, the form with negative part~\eqref{eq:21diag} has a shifted value $0<Q<d$.
\end{lemma}
\begin{proof}
Interchanging $x_4$ and $x_5$ leaves three families. Reflect $x_3$ so that $0\le c:=c_3\le\frac12$. Put $x_2=n$, $x_3=x+c$, $x_4=y+c_y$, $x_5=z+c_z$, $\alpha=a_2$, $b=a_3$, $\kappa=c_1+bc+ec_y+fc_z$, $U=ac$ and $V=ac^2$. Then
\[
Q=nj+\alpha n^2+\kappa n+nbx-ax^2-2Ux-V+n(ey+fz)-\tfrac12\bigl((y+c_y)^2+(z+c_z)^2\bigr),
\]
which is affine in $(d,a,U,V,b,\alpha,\kappa)$. Three Farkas trees of component~I, with witnesses $1\le|n|\le3$, prove the claim on the outer polytope of Appendix~\ref{app:diag21}.
\end{proof}

\begin{theorem}\label{thm:21}
The theorem holds when $(m,K)=(2,1)$; in fact $0<Q<d$ is attained.
\end{theorem}
\begin{proof}
Lemma~\ref{lem:21ternary} covers $A+\delta\le1$. Otherwise the analysis above leads to~\eqref{eq:21diag}, which Lemma~\ref{lem:21diag} contradicts.
\end{proof}

\begin{remark}\label{rem:21transfer}
The resonance method of Section~\ref{sec:22} (exact $t$-partitions, dyadic covers and the ladders of Theorem~\ref{thm:fail}) also applies to much of the low-$A$ region of the case $(2,1)$. For the proof presented here the direct certificate of Lemma~\ref{lem:21ternary} is simpler, and we use it.
\end{remark}
\section{\texorpdfstring{The cases $(1,1)$ and $(1,2)$}{The cases (1,1) and (1,2)}}\label{sec:m1}

When $m=1$ the first completion has width $\delta=d$, and no slack is gained. Apart from the branch $(1,2)$, $L=4$, we therefore work directly with the ternary expansion~\eqref{eq:shiftexp} and with the non-strict goal $0<Q\le d$.

\subsection{\texorpdfstring{The case $(1,1)$}{The case (1,1)}}
\begin{theorem}\label{thm:11}
Suppose $a\ge\frac d2$, $A\ge\frac d2$ and $a+d\le1$. Then for all $\beta\in[-\frac12,\frac12]^3$ and $\nu\in[0,1]$ there are $j\in\Z$ and $w\in\{-3,\dots,3\}^3$ with
\begin{equation}\label{eq:11goal}
0<j+\nu+\beta\cdot w-\phi(w)\le d .
\end{equation}
In particular the theorem holds when $(m,K)=(1,1)$, with $x_2=1$.
\end{theorem}
\begin{proof}
The hypotheses and~\eqref{eq:red} give
\[
\tfrac{\sqrt3}4\le d\le\tfrac23,\qquad t\le2d^3,\qquad A\le\tfrac2{\sqrt3}d^2,\qquad a\le\rho_a<\tfrac{413}{1000},
\]
where $\rho_a^3=(1-\rho_a)^5$; the rational bound follows from $(\frac{413}{1000})^3>(\frac{587}{1000})^5$. For the lower bound on $d$ use $2aAt\ge2\cdot\frac d2\cdot\frac d2\cdot\frac{3d}8$. Define $L$ by $L<d/t\le L+1$. Since $\frac98\le\frac1{2d^2}\le\frac dt\le\frac83$, we have $L\in\{1,2\}$. For $L=1$, $d\ge\frac12$. For $L=2$, $t<\frac d2$, $A\le\frac43t<\frac23d$ and $a\le1-d$ give $3d^3+2d-2<0$, so $d<\rho_2=0.62817\ldots$. We use $I_1=[\frac12,\frac23]$ and $I_2=[\frac{433}{1000},\frac{629}{1000}]$, and for $L=2$ also $A\le\frac{2d}3$ and $a\le\frac{8d}9$. The reflection $x\mapsto-x$ maps $(h,k,\beta_1)$ to $(-h,-k,-\beta_1)$, so we may take $0\le h\le\frac12$ and treat $k\in[0,\frac12]$ and $k\in[-\frac12,0]$ separately. Four Farkas trees of component~I, one for each $(L,\text{sign of }k)$, prove~\eqref{eq:11goal} on the outer polytopes of Appendix~\ref{app:ternary}; they use the strict failure branch $f>d$.
\end{proof}

\begin{remark}
The non-strict inequality in~\eqref{eq:11goal} cannot be improved for the witness family with $x_2=1$. For $\phi=\frac14(w_1^2+w_2^2+w_3^2)$, $d=\frac12$, $\beta=(\frac14,\frac14,\frac14)$ and $\nu=0$, every value $j+\frac14\sum(w_i-w_i^2)$ lies in $\frac12\Z$. This concerns only the restricted family. Other values of $x_2$ can give smaller positive values; for the corresponding unshifted form $Q=x_1x_2+\frac14x_2(x_3+x_4+x_5)-\frac14(x_3^2+x_4^2+x_5^2)$ one has $Q(0,2,1,0,0)=\frac14<d$.
\end{remark}

\subsection{\texorpdfstring{The case $(1,2)$}{The case (1,2)}}
\begin{theorem}\label{thm:12}
The theorem holds when $(m,K)=(1,2)$, i.e.\ when $\frac d2\le a<d$ and $\frac d3\le A<\frac d2$.
\end{theorem}
\begin{proof}
Here $\Delta=d+\frac34A$, and by rounding applied to $H$ we may assume $\Delta+\frac{L^2-1}4t\le1$. The domain satisfies
\[
\tfrac d2\le a<\tfrac{2d}3,\quad \tfrac{3d}8\le A<\tfrac d2,\quad \tfrac1{2\sqrt2}\le d\le\tfrac23,\quad \tfrac32d^3<t\le\tfrac83d^3,\quad A\le\tfrac2{\sqrt3}d^2,\quad a<\tfrac{413}{1000},
\]
and $\frac{41d}{32}\le\Delta<\frac{11d}8$ and $\frac{1107}{1024}\le\frac{\Delta}{t}<5$. Hence $L\in\{1,2,3,4\}$ and $d^2\ge\max(\frac18,\frac{123}{256(L+1)})$. For $L\ge2$, the budget and $t\ge\frac34A$ give $A\le\frac{16(1-d)}{3(L^2+3)}$ and $27(L^2+3)^3d^5\le8192(1-d)^3$. The resulting $d$-intervals are $[\frac{4901}{10^4},\frac23]$, $[\frac{4001}{10^4},\frac{579904}{10^6}]$ and $[\frac{3535}{10^4},\frac{478051}{10^6}]$ for $L=1,2,3$.

\emph{The branch $L=4$.} Since $\Delta\ge\frac{41d}{32}$, $t\ge\frac\Delta5$ and $d\ge\frac1{2\sqrt2}$,
\[
\Delta+t\ge\tfrac65\Delta\ge\tfrac{123}{80}d\ge\frac{123}{160\sqrt2}>\frac12 .
\]
Lemma~\ref{lem:highq} gives $0<G<\Delta$, hence $0<Q<d$ by Lemma~\ref{lem:lift}.

\emph{Ternary regions.} Fourteen Farkas trees of component~I ($15{,}790$ nodes) prove $0<j+\nu+\beta\cdot w-\phi(w)\le d$ on the outer polytopes of Appendix~\ref{app:ternary} intersected with each of the following conditions: $L=2$; $L=1$ with $\Delta+t\le1$; $L=1$, $t=\frac12$, $d+4A\le\frac32$; $L=1$, $t=\frac12$, $d-4A\le-\frac12$; $L=3$ with $\Delta+9t\le2$; $L=3$ with $\Delta-9t\le-1$.

\emph{Remaining part of $L=1$.} Outside these regions $\Delta+t>1$ and $t<\Delta$, so $(\frac12,1)$ applied to $G$ forces $t=\frac12$ and $(a_5'',c_5)\in\{(\frac12,0),(0,\frac12)\}$. Then $|1-4A|+\frac12<d$, and $(1,2)$ applied to $F$ forces $A=\frac14$, $\lambda=0$ and $(a_4',c_4)\in\{(0,0),(\frac14,\frac12)\}$. Now $a=4d^5$, $d\ge8^{-1/4}>\frac7{12}$ and $\frac{29}{100}<a\le\frac13$. As in Section~\ref{sec:22}, $F$ is one of four fixed families, now with $\mu=\frac a4\in[\frac{29}{400},\frac1{12}]$ and target $\frac7{12}<d$. Four dyadic covers give $0<F<\frac7{12}<d$.

\emph{Remaining part of $L=3$.} Here $|\frac32-9t|+\frac12<\Delta$, so failure forces $t=\frac16$ and $(e,c)=(a_5'',c_5)\in\{(\frac16,0),(0,\frac12)\}$. Then $aA=3d^5$, $A\le\frac29$, $a^2\le\frac A3$ and $16^{-1/3}\le d\le(\frac2{27})^{3/10}$. With $q=a_4'/(2A)$ and $C=c_4+\lambda c$, the form $F$ equals~\eqref{eq:corner} with $t_0=\frac16$. Here the linear parameters are treated as free: the coordinates $(\alpha_G-Q_2,\ \kappa_G+ec+2X)$ are replaced by $(\alpha,\kappa)\in[-\frac12,\frac12]^2$. Two Farkas trees of component~I prove $0<F\le d$ on the outer polytope of Appendix~\ref{app:binary12}.
\end{proof}
\section{Proof of Theorem~\ref{thm:main}, sharpness, and verification}\label{sec:end}

\begin{proof}[Proof of Theorem~\ref{thm:main}]
By (I1) and (I2) we may assume the canonical form~\eqref{eq:Q}--\eqref{eq:red}. Proposition~\ref{prop:nonint} treats $c_2\not\equiv0$. For $c_2\equiv0$, Proposition~\ref{prop:largea} treats $a\ge\frac12$, and Proposition~\ref{prop:residual} reduces the remaining forms to the nine pairs of Table~\ref{tab:cases}. These are settled by Theorems~\ref{thm:41}, \ref{thm:31}, \ref{thm:32}, \ref{thm:24}, \ref{thm:23}, \ref{thm:22}, \ref{thm:21}, \ref{thm:11} and~\ref{thm:12}. In every case some $X\equiv c$ satisfies $0<Q(X)\le d=(8|D|)^{1/5}$.

For sharpness, let
\[
Q_*=x_1x_2-\tfrac14\bigl(x_2^2+x_3^2+x_4^2+x_5^2\bigr),\qquad c=(\tfrac12,\dots,\tfrac12),\qquad |D_*|=\tfrac1{256}.
\]
Writing $x_i=u_i/2$ with $u_i$ odd gives $16Q_*=4u_1u_2-\sum_{i=2}^5u_i^2\equiv4-4\equiv0\pmod 8$. So every positive value is at least $\frac12$. The value $\frac12$ is attained at $(\frac32,\frac12,\frac12,\frac12,\frac12)$, and $(\frac12)^5/|D_*|=8$. Hence no constant smaller than $8$ is possible.
\end{proof}

\subsection{The verification data}\label{sec:verify}
All finite data used in the proof are contained in a single release~[\url{https://doi.org/10.5281/zenodo.23789818}], which comprises the two components of Table~\ref{tab:inventory}. A single entry point, \texttt{verify\_release.py}, first checks every file against a manifest of SHA-256 hashes and then runs the exact checkers of both components; per-certificate counts are recorded in the manifest.

\begin{table}[ht]
\centering\small
\begin{tabular}{llrr}
\toprule
Component & Use & Items & Nodes\\
\midrule
I & $(2,1)$: Farkas trees (Lemmas~\ref{lem:21ternary},~\ref{lem:21diag}) & 5 & 1{,}547\\
I & $(1,1)$: Farkas trees (Theorem~\ref{thm:11}) & 4 & 4{,}982\\
I & $(1,2)$: Farkas trees (Theorem~\ref{thm:12}) & 16 & 17{,}224\\
\midrule
II & large $a$: dyadic covers (Proposition~\ref{prop:largea}) & 72 & 5{,}184\\
II & $(3,1)$, $(3,2)$, $(2,4)$, $(2,3)$: dyadic covers & 12 & 1{,}276\\
II & $(2,2)$ and $(1,2)$: dyadic covers & 14 & 1{,}734\\
II & $(2,2)$: Farkas trees at the pinch points (Theorem~\ref{thm:fail}) & 3 & 713\\
II & gap data, $t$-partitions and ladder data (Lemma~\ref{lem:ladder}) & --- & ---\\
\midrule
& Total: Farkas trees / dyadic covers & 28 / 98 & 24{,}466 / 8{,}194\\
\bottomrule
\end{tabular}
\caption{Inventory of the finite data.}\label{tab:inventory}
\end{table}

\emph{Component I} consists of $25$ Farkas trees ($23{,}753$ nodes, $11{,}889$ leaves) with their model generators and exact checkers; the outer polytopes are specified in Appendix~\ref{app:models}. Its checker replays every node and Farkas leaf in exact rational arithmetic.

\emph{Component II} is checked by \texttt{verify\_all.py}, using only the Python standard library; it takes about fifteen seconds. The script first checks that the certificate directory matches the family specification exactly. For each dyadic cover it then reconstructs every box from its path, recomputes the witness constant from the family definition, and verifies the leaf inequalities exactly. It requires four or two children at every internal node and checks restricted targets against their neighbourhoods. For the gap data it regenerates the collision arrangements and all residues from the stated intervals. For the $t$-partitions it checks that the intervals tile the stated ranges, recomputes every lower bound $\delta_-$, every gap maximum and every Macbeath width, and flags exactly the four resonant intervals. It verifies the ladder data of Theorem~\ref{thm:fail}: every rung, the polynomial identities for the infinite rung families, the monotonicity of their ratios, the handoff inequalities and the radii. It checks that every pinch leaf lies in its ball. For the Farkas trees at the pinch points it regenerates the root inequalities from their definition and replays every node and Farkas leaf exactly. Finally it checks the scalar comparisons used in the text. Rational lower bounds for cube roots are seeded by floating-point estimates and then certified by exact cubing. Linear programming was used only to discover the corner trees and plays no role in their verification.

\subsection{Scope}
The finite verifications establish the stated covering properties on the stated rational domains. That every admissible parameter point lies in the appropriate domain is proved in the text and in Appendix~\ref{app:models}. These arguments, together with the classical inputs (I1), (I2) and Lemmas~\ref{lem:dumir}--\ref{lem:jackson}, are not formalised in a proof assistant. The finite certificates are elementary to check, and an independent reimplementation of the checkers should need only the definitions of Section~\ref{subsec:certs} and the models of Appendix~\ref{app:models}. The analytic inclusions are the natural target for a formal verification of Theorem~\ref{thm:main}.
\section{Concluding remarks}\label{sec:concl}

\emph{What is settled.} Theorem~\ref{thm:main} determines the last of the constants $\Gamma_{r,s}$ for nondegenerate indefinite forms. Its proof uses the classical canonical reduction and the gap lemmas of Dumir, Macbeath and Jackson, as earlier work did. What is new is the treatment of the subcases where those lemmas leave genuine exceptional configurations; this concerns the four pairs $(m,K)$ in the last three rows of Table~\ref{tab:status}, where the earlier bounds were weaker than~$8$.

\emph{Equality.} Section~\ref{sec:end} shows sharpness by one explicit form and shift. We have not classified the pairs $(Q,c)$ for which no $x\equiv c$ satisfies $0<Q(x)<(8|D|)^{1/5}$. Bhardwaj, Kathuria and Raka exhibit critical forms in several of the subcases they treat~\cite{BKR}, and isolation results of this kind are known for other types~\cite{BDH84}. A complete description of the critical pairs for type $(1,4)$, and whether $Q_*$ is isolated among them, remains open.

\emph{Methods that transfer.} Several ingredients are independent of the type $(1,4)$.
\begin{itemize}
\item The ladder lemma and the completeness criterion (Lemma~\ref{lem:ladder}, Proposition~\ref{prop:complete}) separate two phenomena that a bounded search cannot distinguish: regions where successful witnesses cannot be chosen uniformly bounded as the parameters approach an exceptional point, and points where the reduced family genuinely fails. Once the exact failure set is known, one can return to a less reduced form at those points only. We expect this pattern, certifying the failure set of the reduced problem and then restoring the discarded coupling there, to be useful in other inhomogeneous problems settled by successive square completion, for instance second minima or asymmetric inequalities.
\item The quotient-pruning lemma (Lemma~\ref{lem:prune}) holds in every dimension. Remark~\ref{rem:six} records that for type $(1,5)$ it leaves five pairs $(m,K)$, one of which is closed by the argument of Section~\ref{sec:K1}. Since $\Gamma_{1,5}$ is known, this is only a test of the method; whether the remaining four pairs yield to the tools of Sections~\ref{sec:third}--\ref{sec:22} would show how far the approach reaches in higher dimension.
\item Proposition~\ref{prop:gapexact} computes uniform residue gaps exactly. It applies wherever one-variable gap arguments of Macbeath's type are used over an interval of coefficients, and it replaces case-by-case gap estimates by a finite line arrangement.
\end{itemize}

The useful distinction is therefore between failure of the original inequality and failure introduced by reduction: identifying the latter exactly makes it possible to recover the missing coupling only where it is needed.
\appendix
\section{Outer polytopes for the Farkas trees}\label{app:models}

Every Farkas tree in this paper is proved on a rational polytope that contains all admissible parameter points. We list the inequalities. Each follows from the reduction inequalities, from the case hypotheses, or from one of the following elementary facts.
\begin{itemize}
\item[(M1)] \emph{Interval moments.} If $l\le x\le u$ and $\omega>0$, then $P=\omega x$ and $R=\omega x^2$ satisfy $l\omega\le P\le u\omega$, $R\le(l+u)P-lu\omega$ and $R\ge2\xi P-\xi^2\omega$ for every $\xi$.
\item[(M2)] \emph{Products.} If moreover $r_0\le y\le s_0$, $T=\omega y$ and $S=\omega xy$, then
$S\ge r_0P+lT-lr_0\omega$, $S\ge s_0P+uT-us_0\omega$, $S\le r_0P+uT-ur_0\omega$ and $S\le s_0P+lT-ls_0\omega$.
\item[(M3)] \emph{Squares of combinations.} $\omega(x-\xi y)^2\ge0$ gives $R_x-2\xi S+\xi^2R_y\ge0$.
\item[(M4)] \emph{Convex functions of $d$.} If $d$ ranges over $[\ell,H]$ and $g$ is convex, then $g$ lies below its secant and above its tangents. For $g(d)=d^{5/3}$ the secant endpoints are rounded upward to multiples of $10^{-6}$ by exact cubing; $U(q)$ denotes the least multiple of $10^{-6}$ whose cube is at least $q$.
\item[(M5)] \emph{Determinant cuts.} On a coefficient box for $(a,A,t)$, every affine function interpolating $aAt$ at four independent corners that bounds $aAt$ from below (above) at all eight corners bounds it on the whole box, by multiaffinity. With $S_5$ the secant of $d^5$ on $[\ell,H]$, we impose $2P_-\le S_5(d)$ and $2P_+\ge5\xi^4d-4\xi^5$ for $\xi=\ell,\frac{3\ell+H}4,\frac{\ell+H}2,\frac{\ell+3H}4,H$.
\end{itemize}

\subsection{The ternary model}\label{app:ternary}
For the expansion~\eqref{eq:shiftexp} put
\[
(p,q,r,s,u,b,v)=(ah,\,ak,\,ah^2,\,ahk,\,ak^2,\,A\lambda,\,A\lambda^2).
\] In the fifteen coordinates $\theta=(d,a,A,t,p,q,r,s,u,b,v,\beta_1,\beta_2,\beta_3,\nu)$ the Gram matrix of $\phi$ is
\[
\begin{pmatrix}a&p&q\\p&A+r&b+s\\q&b+s&t+u+v\end{pmatrix},
\]
and every witness $j+\nu+\beta\cdot w-\phi(w)$ is affine in $\theta$. The common inequalities are: $a+d\le1$, $A\ge\frac34a$, $t\ge\frac34A$, $A\le t+v$, $\beta\in[-\frac12,\frac12]^3$ and $\nu\in[0,1]$. Further, $|p|,|q|\le\frac a2$, $r,u\le\frac a4$, $r\ge2\xi p-\xi^2a$ and $u\ge2\xi q-\xi^2a$ for $\xi\in\{0,\pm\frac14,\pm\frac12\}$, and $-\frac a4+\frac{\epsilon p+\eta q}2-\epsilon\eta s\le0$ for $\epsilon,\eta=\pm1$. Also $0\le b\le\frac A2$, $v\le\frac b2$ and $v\ge2\xi b-\xi^2A$ for $\xi\in\{0,\frac18,\frac14,\frac38,\frac12\}$. Finally $\phi(w)\ge a$ for one vector $w$ from each pair $\pm w$ with $0<\|w\|^2\le5$. The shear intervals are $h\in[0,\frac12]$ (after reflection) and $k\in[0,\frac12]$ or $[-\frac12,0]$, and (M1)--(M2) are applied on them. The case-specific inequalities are as follows.

\emph{$(2,1)$, low $A$.} The hypotheses of Section~\ref{sec:21}, \eqref{eq:21bounds} and $A+\delta\le1$, with $\ell=\frac{4419}{10^4}$, $H=\frac8{15}$ and $R=\frac{1414214}{10^6}>\sqrt2$. Secants of $a\le d^{5/3}$, $A\le Rd^2$ and $t\le\frac{12}5d^3$ on $[\ell,H]$. Determinant cuts on the box $a\in[\frac\ell3,\min(\frac H2,1-\ell,U(H^5))]$, $A\in[\frac{5\ell}8,\min(\frac{11H}8,RH^2,1-\frac{5\ell}4)]$, $t\in[\frac{15\ell}{32},\frac{12}5H^3]$.

\emph{$(1,1)$.} The hypotheses of Theorem~\ref{thm:11}: $d\in I_L$, $\frac d2\le a\le d$, $\frac d2\le A\le d$, $a\le\frac{413}{1000}$, $Lt\le d\le(L+1)t$, and for $L=2$ also $A\le\frac{2d}3$ and $a\le\frac{8d}9$. With $I_L=[\ell,H]$ and $R=\frac{577351}{500000}>\frac2{\sqrt3}$: secants of $t\le2d^3$, $A\le Rd^2$ and $a\le d^{5/3}$. Determinant cuts on $a\in[\frac\ell2,\min(H,1-\ell,\frac{413}{1000},U(H^5))]$, $A\in[\frac\ell2,\min(H,RH^2,B_L)]$ with $B_1=H$, $B_2=\frac{2H}3$, and $t\in[\max(\frac{3\ell}8,\frac\ell{L+1}),\min(2H^3,\frac HL)]$. The trees use the non-strict goal, i.e.\ the strict failure row $f>d$.

\emph{$(1,2)$.} The hypotheses and bounds of Theorem~\ref{thm:12} on the stated $d$-interval $[\ell,H]$, together with $L t\le\Delta\le(L+1)t$ and the budget. For $L\ge2$ add also $d\le\frac{128}{41(L+3)}$, $d\le\frac{128}{155+9L^2}$, $3(L-1)A\le4d$ and $9(L-1)a\le16d$. Secants of $t\le\frac83d^3$, $A\le Rd^2$ and $a\le d^{5/3}$, and tangents of $t\ge\frac32d^3$ at $\ell,\frac{\ell+H}2,H$. Determinant cuts on $a\in[\frac\ell2,\min(\frac{2H}3,1-\ell,\frac{413}{1000},U(H^5))]$, $A\in[\frac{3\ell}8,\min(\frac H2,RH^2)]$ and
\[
t\in\Bigl[\max\bigl(\tfrac{9\ell}{32},\tfrac{41\ell}{32(L+1)},\tfrac32\ell^3\bigr),\ \min\bigl(\tfrac83H^3,\tfrac{11H}{8L},T_L\bigr)\Bigr],
\]
with $T_1=\frac{11H}8$ and $T_L=\frac{4(1-41\ell/32)}{L^2-1}$ for $L\ge2$. In the region $L=3$, $\Delta+9t\le2$, the interval of $\lambda$ is split into $[0,\frac14]$ and $[\frac14,\frac12]$.

\subsection{\texorpdfstring{The $(2,1)$ diagonal model}{The (2,1) diagonal model}}\label{app:diag21}
Coordinates $(d,a,U,V,b,\alpha,\kappa)$ with $\frac d3\le a\le\frac d2$, $a+d\le1$, $\frac{6389}{10^4}\le d\le\frac7{10}$ (as $6^{-1/4}>0.6389$ and $2(\frac7{10})^5+\frac7{10}>1$), and $a\le\frac5{16}$ (if $a\ge\frac5{16}$ then $d\le\frac{11}{16}$ and $2(\frac{11}{16})^5<\frac5{16}$). Further $0\le U\le\frac a2$, $V\le\frac U2$, $V\ge2\xi U-\xi^2a$ for $\xi\in\{0,\frac18,\frac14,\frac38,\frac12\}$, and $b,\alpha,\kappa\in[-\frac12,\frac12]$. The relation $a=2d^5$ enters through its secant and nine equally spaced tangents.

\subsection{\texorpdfstring{The $(1,2)$ coupled model}{The (1,2) coupled model}}\label{app:binary12}
Coordinates $(d,a,A,p,v,r,s,u,w,m,\alpha,\kappa)$ as in~\eqref{eq:corner}, with $t_0=\frac16$ and free $(\alpha,\kappa)\in[-\frac12,\frac12]^2$. Constraints: $d\in[\frac{3968}{10^4},\frac{458037}{10^6}]$, $\frac d2\le a\le\frac{2d}3$, $\frac{3d}8\le A\le\frac d2$, $a+d\le1$, $A\ge\frac34a$, $A\le\frac29$, $a\le U_2(\frac2{27})$, $\frac12\le d+\frac34A\le\frac23$ and $A-v\le\frac16$. Apply (M1)--(M3) to $\lambda\in[0,\frac12]$ and $q,C\in[-\frac12,\frac12]$, with five square tangents and $m-2\xi w+\xi^2v\ge0$ for $\xi\in\{0,\pm\frac12,\pm1,\pm2\}$. Add the secant of $a\le d^{5/3}$, McCormick bounds for $aA$ on $a\in[\frac\ell2,\min(\frac{2H}3,U_2(\frac2{27}),U_3(H^5))]$ and $A\in[\frac{3\ell}8,\min(\frac H2,\frac29)]$, and $aA=3d^5$ through the secant and five tangents of $d^5$. Here $U_j(x)$ rounds $x^{1/j}$ upward to a multiple of $10^{-6}$.

\subsection{\texorpdfstring{The $(2,2)$ corner models}{The (2,2) corner models}}\label{app:corner}
Fix $(L,t_0,e,c)$ and a pinch point $(\alpha^*,\kappa^*)$ from Theorem~\ref{thm:fail}. Let $d_l$ be the rational lower bound for $(\frac{5t_0}{18})^{1/3}$ used in Section~\ref{sec:22}, $d_h=\frac47\tau$ and $\tau=\frac34$. The coordinates are those of~\eqref{eq:corner}, and the inequalities are as follows.
\begin{enumerate}
\item $d\in[d_l,d_h]$. The upper bound holds because $d+a+A\ge\frac74d$ on the $(2,2)$ domain, and the lower bound because $t_0=d^5/(2aA)\le\frac{18}5d^3$.
\item $\frac d3\le a\le\frac d2$, $\frac13(d+\frac34a)\le A\le\frac12(d+\frac34a)$, $A\ge\frac34a$ and $\frac34A\le t_0$.
\item $Lt_0\le d+\frac34(a+A)\le\min\bigl((L+1)t_0,\,1-\frac{L^2-1}4t_0\bigr)$, and $d+a+A\le\tau$. For $L=3$ also $\frac{a+A}4\ge\frac1{12}$.
\item $P=aA$: the four McCormick bounds on $a\in[\frac{d_l}3,\frac{d_h}2]$, $A\in[\frac{5d_l}{12},\min(\frac43t_0,\tau-\frac43d_l)]$. Also $P\le$ the secant of $d^5/(2t_0)$ on $[d_l,d_h]$, and $P\ge$ its tangents at seven equally spaced points.
\item $\lambda\in[0,\frac12]$: $0\le p\le\frac A2$, $v\le\frac p2$, $v\ge2\xi p-\xi^2A$ for $\xi\in\{0,\frac18,\frac14,\frac38,\frac12\}$, and $v\ge A-t_0$ (from $t_0+A\lambda^2\ge A$).
\item $q,C\in[-\frac12,\frac12]$: $|r|,|u|\le\frac A2$; (M2) for $s$ (with $p,r$) and $w$ (with $p,u$); $Q_2\le\frac A4$, $Q_2\ge0$ and $Q_2\ge2\xi r-\xi^2A$, and $m\le\frac A4$ and $m\ge2\xi u-\xi^2A$, for $\xi\in\{0,\pm\frac14,\pm\frac12\}$; (M2) for $X$ (with $r,u$).
\item (M3): $m-2\xi w+\xi^2v\ge0$, $v-2\xi s+\xi^2Q_2\ge0$ and $m-2\xi X+\xi^2Q_2\ge0$ for $\xi\in\{0,\pm\frac12,\pm1,\pm2\}$.
\item $\alpha_G=\alpha^*$ and $\kappa_G=\kappa^*$.
\end{enumerate}
The target is $T=d+\frac34a$, with closed failure rows $f\le0$ and $f\ge T$. The generator of these rows is part of component~II, and the checker regenerates them before replaying the trees.

\end{document}